\documentclass[10pt,twoside,a4paper]{article}

\usepackage[utf8]{inputenc}
\usepackage[T1]{fontenc}
\usepackage[english]{babel}

\usepackage{amsmath}
\usepackage{amsfonts}
\usepackage{amssymb}
\usepackage{amsthm}
\usepackage{amsxtra}
\usepackage{mathtools}
\usepackage[stddef=false]{delimset}

\usepackage{authblk}
\usepackage[title,titletoc]{appendix}
\usepackage{csquotes}
\usepackage{url}
\usepackage{enumitem}
\usepackage{environ}

\usepackage{graphicx}
\graphicspath{{./images/}}
\usepackage{subcaption}

\usepackage[margin=2cm]{geometry}

\usepackage[
    backend=biber,
    style=numeric,
    giveninits=true,
    maxbibnames=99,
    sorting=nyt,
    url=false,
    doi=true,
    eprint=true
]{biblatex}
\usepackage{hyperref}

\makeatletter
\newcommand\blfootnote[1]{%
  \begin{NoHyper}%
  \renewcommand\thefootnote{}\footnote{#1}%
  \addtocounter{footnote}{-1}%
  \end{NoHyper}%
}

\newcommand{\email}[1]{\@ifundefined{@emails}{\gdef\@emails{\href{mailto:#1}{#1}}}{\g@addto@macro\@emails{, \href{mailto:#1}{#1}}}}
\newcommand{\emails}[1]{\def\@emails{#1}}
\newcommand{\printemails}{\@ifundefined{@emails}{}{\noindent\textbf{Email addresses:} \@emails}}

\newcommand{\keywords}[1]{\def\@keywords{#1}}
\newcommand{\printkeywords}{\@ifundefined{@keywords}{}{\noindent\textbf{Keywords and phrases:} \@keywords.}}

\newcommand{\subjclass}[1]{\def\@subjclass{#1}}
\newcommand{\printsubjclass}{\@ifundefined{@subjclass}{}{\noindent\textbf{Mathematics Subject Classification (2020):} \@subjclass.}}
\makeatother

\numberwithin{equation}{section}
\numberwithin{figure}{section}

\newtheoremstyle{thm-style-oskari}{7pt}{7pt}{\itshape}{}{\scshape}{.}{.5em}{}
\theoremstyle{thm-style-oskari}
\newtheorem{theorem}{Theorem}[section]

\newtheorem{corollary}[theorem]{Corollary}
\newtheorem{lemma}[theorem]{Lemma}

\newtheorem{remark}[theorem]{Remark}

\NewEnviron{eqsplit}{%
    \begin{equation} \begin{split}
        \BODY
    \end{split} \end{equation}
}

\NewEnviron{eqsplit*}{%
    \begin{equation*} \begin{split}
        \BODY
    \end{split} \end{equation*}
}

\newcommand{\ol}[1]{\overline{#1} \!\,}
\newcommand{\wt}{\widetilde}
\newcommand{\eps}{\varepsilon}

\newcommand{\E}{\mathbb{E}}
\newcommand{\R}{\mathbb{R}}
\newcommand{\C}{\mathbb{C}}
\newcommand{\N}{\mathbb{N}}
\newcommand{\Z}{\mathbb{Z}}
\newcommand{\ee}{\mathrm{e}}
\newcommand{\ii}{\mathrm{i}}
\newcommand{\dd}{\mathrm{d}}

\newcommand{\p}{\delim()}
\newcommand{\s}{\delim[]}
\newcommand{\abs}{\delim{\lvert}{\rvert}}
\newcommand{\ceil}{\delim{\lceil}{\rceil}}

\DeclareMathOperator{\diag}{diag}
\DeclareMathOperator{\Id}{\mathbb{I}}
\DeclareMathOperator{\Tr}{Tr}
\DeclareMathOperator{\Str}{Str}
\let\Re\relax\DeclareMathOperator{\Re}{Re}
\let\Im\relax\DeclareMathOperator{\Im}{Im}

\DeclareMathOperator*{\spec}{Spec}
\let\det\relax\DeclareMathOperator{\det}{det}
\DeclareMathOperator{\Sdet}{Sdet}

\newcommand{\1}{\mspace{1 mu}}

\newcommand{\HypergeometricPFQ}[5]{
    {}_{#1} F_{#2}\delimpair{(}{[m]|}{)}*{
        \begin{matrix} #3 \\ #4 \end{matrix}
    }{
        #5
    }
}

\newcommand{\MeijerG}[5]{
    G^{#1}_{#2}\delimpair{(}{[m]|}{)}*{
        \begin{matrix} #3 \\ #4 \end{matrix}
    }{
        #5
    }
}

\newcommand*{\DLMFsecref}[1]{\href{https://dlmf.nist.gov/#1}{\S #1}}
\newcommand*{\DLMFeqref}[2]{\href{https://dlmf.nist.gov/#1.E#2}{(#1.#2)}}

\title{Universal scaling limits at the spectral singularity of structured random matrices}
\author[1]{Markus Ebke}
\email{markus.ebke@fau.de}
\author[1]{Torben Krüger}
\email{torben.krueger@fau.de}
\affil[1]{Department of Mathematics, Friedrich-Alexander-Universität Erlangen-Nürnberg, Cauerstraße 11, 91058 Erlangen, Germany}

\keywords{Wigner-type matrix; eigenvalue distribution; microscopic scaling limit; supersymmetric method; superbosonization}
\subjclass{%
Primary: 60B20; 
60F15; 
Secondary: 33C20; 
33C60
}

\makeatletter
\hypersetup{
    pdftitle={\@title},
    pdfauthor={Markus Ebke, Torben Krüger},
    pdfsubject={\@subjclass},
    pdfkeywords={\@keywords},
}
\makeatother

\begin{document}
\date{}
\maketitle
\thispagestyle{empty}

\begin{abstract}
The empirical spectral distribution of Hermitian $K \times K$-block random matrices converges to a deterministic density on the real line with a potential atom at the origin as the dimension of the blocks tends to infinity. In this model the variances of the entries depends on the block and the limiting density is determined by these variances. In the absence of an atom the density is either bounded or has a power law singularity at the origin. We determine all scaling limits of the spectral density on the eigenvalue spacing scale at this singularity for Gaussian matrices with block numbers $K \leq 3$. The appropriate scaling for the universal limit is correctly predicted by the global eigenvalue density. For $K=3$ the local one-point function exhibits an additional logarithmic singularity. The scaling limits depend only on the zero pattern within the variance profile, but not on the values of its positive entries.
\end{abstract}

\printkeywords\par
\printsubjclass
\blfootnote{\printemails}


\section{Introduction}
A prominent theme within the theory of random matrices (RMT) is the phenomenon of universality: the empirical spectral distribution of large random matrices depends only on a few model features, such as the underlying symmetry class of the matrix or the presence of global constraints on the entries, not however on the details of the entry distributions. This phenomenon is observed on two complementary scales. On the global (macroscopic) scale, the empirical eigenvalue distributions converge to deterministic laws that reflect coarse structural information about the ensemble such as the covariance structure of the entries. Wigner’s semicircle law \cite{wigner1958DistributionRootsSymmetric} and the Marchenko–Pastur law \cite{marcenko1967DistributionEigenvaluesSets} are classical examples for which the limit of the empirical spectral measure depends only on the entry variances but not on their distributions. On the local (microscopic) spectral scale, after zooming to the typical eigenvalue spacing, universal correlation functions emerge: the sine kernel in the bulk \cite{mehta1960StatisticalPropertiesLevelspacings,gaudin1961LoiLimiteLespacement}, Tracy–Widom/Airy statistics at soft edges \cite{tracy1994Airy,tracy1996OrthogonalSymplectic}, Bessel-type laws at hard edges \cite{tracy1994Bessel}, and the Pearcey-type correlations at cusps \cite{brezin1998ClosureGap}. Both, global laws and local spectral universality, have by now been established for wide classes of random matrix ensembles far beyond the classical settings of Wigner matrices with i.i.d.\ entries, for which the sine kernel universality was proved in \cite{soshnikov1999UniversalityEdgeWigner,tao2011RandomMatricesUniversality,erdos2012BulkUniversalityGeneralized}.

Within this landscape of random matrix ensembles, \emph{structured} or \emph{Wigner-type} models play a central role, both to test theoretical hypotheses and to model applications with inhomogeneous settings, where the distribution of the entries depends on their position inside the matrix. These ensembles are Hermitian matrices $H=(h_{\alpha \beta})$ with independent entries and a variance profile $\E\1 \abs{h_{\alpha \beta}}^2$. The simplest instance are $N \times N$ Wigner matrices with i.i.d.\ entries whose variances equal $\frac{1}{N}$.
Block variance profiles are of particular interest due to the transparent nature of their limiting spectral measure. In these models the matrix indices are grouped into $K$ batches $I_1, \dots, I_K$ of size $N$ and the variance profile only depends on the batch index, i.e.\ $\E\1 \abs{h_{\alpha \beta}}^2 = \frac{1}{N} s_{ij}$ for all $(\alpha,\beta) \in I_i\times I_j$. As $N$ tends to infinity, the empirical spectral distribution $\rho$ is governed by the associated vector Dyson equation and depends on $S=(s_{ij})$ (see e.g. \cite{ajanki2017SingularitiesSolutions} for details and references).
This framework includes block models and multi-type networks (e.g.\ inhomogeneous Erd{\H o}s–R\'enyi and stochastic block models \cite{lei2015ConsistencySpectralClustering}), block band matrices with applications in condensed matter physics \cite{erdos2013DelocalizationDiffusion}, graph Laplacians with community structure \cite{rohe2011SpectralClusteringBlockmodel,chatterjee2022SpectralPropertiesLaplacian}, bipartite/chiral Hamiltonians \cite{altland1997NonstandardSymmetryClasses}, as well as heteroskedastic noise models in high-dimensional statistics and signal processing \cite{hachem2007DeterministicFunctionals,couillet2011WirelessCommunications,alt2017LocalLawGram}. Systematic analysis of the vector Dyson equation as well as bulk local spectral universality for Wigner-type matrices has been developed in \cite{ajanki2017SingularitiesSolutions,ajanki2017GeneralWigner}.

It is natural to ask for a classification of all local universality classes that can occur within such structured ensembles. We restrict our discussion here to block models with centred entries. Under sufficiently strong non-degeneracy assumptions on the variance profile $S$, e.g.\ when $S$ is strictly positive, it was established in \cite{ajanki2017GeneralWigner} that the limiting spectral measure is bounded, analytic, whenever positive, and approaches zero only at regular square root edges or cubic root cusps. The corresponding universality classes are the sine \cite{ajanki2017GeneralWigner}, Airy \cite{alt2020CorrelatedRandomMatrices} and Pearcey statistics \cite{cipolloni2019CuspUniversalityReal,erdos2020CuspUniversalityComplex}, respectively.
However, once zero blocks within the matrix are permitted $\rho$ may have an atom at the origin or a power law singularity of the form $\rho(E)\sim \abs{E}^{-\sigma}$ for $E \to 0$. These singularities have been classified in \cite{krueger2025SingularityDegree} depending on the zero pattern of the variance profile. In particular, $\sigma = (\ell-1)/(\ell+1)$ for some integer $1\le \ell\le K$. Some special cases are also considered in \cite{Kolupaiev2021AnomalousSingularity}, which additionally considers $N$-dependent variance profiles $S$ with entries that are uniformly bounded from above and below.
A related model is the linearization of covariance matrices $X X^*$, where $X$ are $\lambda$-shaped matrices, which are rectangular with independent complex entries whose variance profile has ones and zeros in the shape of a Young-diagram. The study of the eigenvalue distribution of these ensembles is closely related to the combinatorics of plane trees \cite{cunden2023RandomMatricesYoung,bisi2025LshapedRandomMatrices}.
The question which local spectral universality class is associated to the singularity with given $\ell$ however has not been addressed in previous works. Only numerical evidence for the conjectured spacing scale $\eta_N \sim N^{-(\ell+1)/2}$ at the singularities for $K=2$ and $K=3$ blocks has been provided in \cite{krueger2025SingularityDegree}.

In the present work we address the question of local spectral statistics for $K=2,3$ with Gaussian entries and derive a formula for the spectral density at the singularity on the eigenvalue spacing scale. For $K=2$ we show that the global power law, dictated by the density $\rho$, persists down to the microscopic scale. With $\ell=1$ this yields a new universality class whose limiting one-point function is expressed in terms of Meijer $G$ and generalized hypergeometric functions.
In the $K=3$ case with $\ell=2$ a new universality class is established at a typical spacing scale $\sim N^{-2}$, whose local one-point function is expressed in terms of modified Bessel functions. In this case, the local density has a logarithmic singularity at the origin suggesting a particularly strong affinity of the smallest singular value to the origin. In both cases the limits depend only on the divergence type of the macroscopic density and not on the specific values of the non-zero entries within the variance profile $S$.
Furthermore, we consider variants of $K=2, 3$ where entries of the variance profile vanish as $N \to \infty$.
This allows us to interpolate between our new results for the $K=2$ case and the well-known results from \cite{forrester1993SpectrumEdge} for the chiral Gaussian Unitary Ensemble.
We also examine an analogous limit for $K=3$.

For the proofs we employ the supersymmetric method, including superbosonization, to derive finite-$N$ representations for the expected trace of the resolvent of $H$.
The supersymmetric method has been used for a variety of random matrix ensembles including random band matrices \cite{disertori2002DensityStatesRandom,disertori2017DensityStatesRandom,shcherbina2014UniversalityLocalRegime,spencer2012SUSYStatisticalMechanics}, interpolating ensembles \cite{shamis2013DensityStatesGaussian,afanasiev2022CorrelationFunctionsCharacteristic}, the shifted Ginibre ensemble \cite{cipolloni2020OptimalLower,cipolloni2022ConditionNumberShifted,shcherbina2022LeastSingularValue,maltsev2025BulkUniversalityComplex} and chiral random matrices \cite{kaymak2014SupersymmetryChiral}.
The superbosonization formula was used in \cite{bunder2007SuperbosonizationFormula} to derive the eigenvalue distribution for almost diagonal matrices, i.e.\ $N = 1$ and $S_{ij} \to 0$ for $\abs{i - j} \to \infty$ in our notation.
To determine the local scaling limit of the eigenvalue density at the singularity of the density of states, we rescale the integrals of the finite-$N$ representations to accommodate for the divergent behaviour of the resolvent at the singularity. Here, additional symmetries in the resulting action at the spacing scale emerge. After a higher order saddle point analysis these symmetries are integrated out leading to the special functions within the final formulas.

The paper is organized as follows.
In Section~\ref{sec: Main Results} we define the matrix model for $H$ and state our main results along with some remarks and plots of the universal local spectral densities.
The following sections contain the proofs:
First, in Section~\ref{sec: Integral Representation} we derive an exact integral representation for the expectation of the resolvent $(H - z)^{-1}$ at finite matrix size.
In Section~\ref{sec: Large N Limit at the Origin} we compute the $N \to \infty$, $z \to 0$ asymptotics of this integral representation for $K = 2$ and $K = 3$.
In Section~\ref{sec: Microscopic Scaling Limit} we use the asymptotic formula to determine the microscopic scaling limit of the one-point function.
We investigate the weak non-chirality regime in Section~\ref{sec: Weak Non-Chirality Limit}, here we start again with the finite $N$ integral representation from the end of Section~\ref{sec: Large N Limit at the Origin} and compute its asymptotics for an $N$-dependent variance matrix $S$.
In Appendix~\ref{app: Saddle-Point Method} we summarize the saddle-point method and state the formulas that we use in Section~\ref{sec: Large N Limit at the Origin} and Section~\ref{sec: Weak Non-Chirality Limit}.

\paragraph{Acknowledgements}
We thank Gernot Akemann and the participants of the PEPPER autumn-school at the University of Münster in November 2024 for enlightening discussions.

\section{Main Results}
\label{sec: Main Results}

We consider Hermitian $K \times K$-block random matrices with block dimension $N$ that have independent centred Gaussian entries.
The block variance profile is $S = (s_{ij})_{i,j=1}^K$ with $s_{ij} = s_{ji} \geq 0$.
This model is conveniently written in tensor notation as
\begin{equation} \label{eq: matrix model definition}
    H = \sum_{i,j=1}^K \sqrt{\frac{s_{ij}}{2}} \, \p*{E_{ij} \otimes X_{ij} + E_{ji} \otimes X_{ij}^*}
    \in \C^{K \times K} \otimes \C^{N \times N} \simeq \C^{KN \times KN}\,,
\end{equation}
where $E_{i j} = (\delta_{i k} \delta_{j l})_{k,l=1}^K$ is the canonical basis of $\C^{K \times K}$ and $X_{ij} \in \C^{N \times N}$ are independent complex $N \times N$ Ginibre matrices.
A Ginibre matrix $X = (x_{kl})_{k, l = 1}^N$ has centred i.i.d.\ Gaussian entries with $\E\1 \abs{x_{kl}}^2 = \frac{1}{N}$ and $\E\1 x_{kl}^2 = 0$.

We define the normalized one-point function (also called the density of states) $\rho_N$ of $H$ and the asymptotic density of states $\rho_\infty$ as
\[
    \rho_N := \frac{1}{K N}\,\E\s*{\sum_{\lambda \in \spec(H)} \delta_{\lambda}}\,, \qquad
    \rho_{\infty} := \lim_{N \to \infty} \rho_N\,,
\]
respectively, where the sum over the eigenvalues $\lambda \in \spec(H)$ of $H$ is taken with multiplicity.
The density of states is connected to the resolvent of $H$ through the Stieltjes transform
\[
    \frac{1}{K N} \E\s*{\Tr\p*{\frac{1}{H - z}}}
    = \frac{1}{K N}\E\s*{\sum_{\lambda \in \spec(H)} \frac{1}{\lambda - z}}
    = \int_{\R} \frac{\rho_N(\dd E)}{E-z}\,,\qquad \text{for } z \in \C \setminus \R.
\]
Here $\frac{1}{H - z}$ denotes the resolvent of $H$, i.e.\ the matrix inverse of $H - z \Id_{K N}$ where $\Id_{K N}$ is the $KN \times KN$ identity matrix.
The trace $\Tr$ is the sum over all diagonal entries.
The density can be recovered using the inversion formula
\begin{equation} \label{eq: Stieltjes inversion}
    \rho_N(E) = \frac{1}{\pi} \lim_{\eps \searrow 0} \frac{1}{K N} \Im \E\s*{\Tr\p*{\frac{1}{H - (E + \ii \eps)}}}
\end{equation}
which holds if the limit on the right-hand side exists.
As a first step towards analyzing the density $\rho_N$, we derive an exact expression for the expectation of the resolvent at finite $N$.
The proof is given in Section~\ref{sec: Integral Representation}.

\begin{lemma} \label{thm: Finite N Integral Representation}
    Let $H$ be given by \eqref{eq: matrix model definition}.
    Then for all $z \in \C$ with $\Im(z) > 0$ and all $N \in \N$ the expectation of the resolvent has the integral representation
    \begin{equation} \label{eq: resolvent expectation integral}
        \E\s*{\Tr\p*{\frac{1}{H - z}}}
        = \frac{\ii N^{K+1}}{(2 \pi \ii)^K} \int_{(0, \infty)^K} \dd a \oint^K \dd b \, \ee^{-N I(a) + N I(b)} \det\s*{ S + \diag_{i = 1}^K\p*{\frac{1}{a_i b_i}} } \sum_{j = 1}^K a_j,
    \end{equation}
    where $\dd a$ is the $K$-dimensional Lebesgue-measure and the $b$-integral is $\oint^K \dd b := \oint \dd b_1 \dots \oint \dd b_K$ where each $\oint \dd b_i$ is a counter-clockwise contour integral over the circle around $0$.
    The function $I: \R^K \to \C$ in the exponent is
    \begin{equation} \label{eq: I definition}
        I(a) := \frac{1}{2} \sum_{i, j = 1}^K s_{ij} a_j a_i - \ii z \sum_{i = 1}^K a_i - \sum_{i = 1}^K \log(a_i).
    \end{equation}
\end{lemma}

\begin{figure}[ht]
    \centering
    \begin{subfigure}[t]{0.4\textwidth}
        \centering
        \includegraphics[height=130pt]{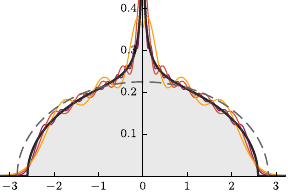}
        \caption{$K = 2$ and $S = \begin{pmatrix}1 & 1 \\ 1 & 0\end{pmatrix}$}
    \end{subfigure}
    \hfill
    \begin{subfigure}[t]{0.59\textwidth}
        \centering
        \includegraphics[height=130pt]{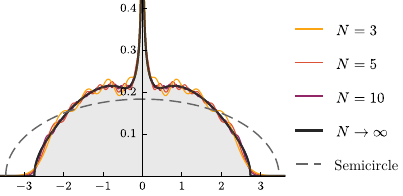}
        \caption{$K = 3$ and $S = \begin{pmatrix}0 & 1 & 1 \\ 1 & 1 & 0 \\ 1 & 0 & 0\end{pmatrix}$}
    \end{subfigure}
    \caption{Plots of the eigenvalue density of the random matrix $H$. The finite $N$ densities (colored lines) are histograms obtained from numerical simulations of \eqref{eq: matrix model definition}. The asymptotic density $\rho_\infty$ (black line) is derived from \eqref{eq: Asymptotic density from saddle points}. For comparison the gray dashed line is the semicircle law that would hold if all entries of $S$ were 1.}
    \label{fig: macroscopic density}
\end{figure}

For Wigner-type matrices it is well known that (under certain assumptions on the entry distribution, see \cite{girko2001TheoryStochastic, ajanki2017SingularitiesSolutions, ajanki2019QuadraticVectorEquations}) the asymptotic density of states $\rho_{\infty}$ can be expressed as
\begin{equation} \label{eq: Asymptotic density from saddle points}
    \rho_{\infty}(E) = \frac{1}{\pi} \lim_{\eps \searrow 0} \Re\frac{1}{K} \sum_{i=1}^K a_i(E + \ii \eps),
\end{equation}
where $a(z) = (a_1(z), \dots, a_K(z))$ is the unique solution of the quadratic vector equation $\nabla I(a(z)) = 0$ such that $\Re a_i(z) > 0$ for all $i = 1, \dots, K$.
In Figure~\ref{fig: macroscopic density} we plot the eigenvalue density of $H$ from \eqref{eq: matrix model definition} for small $N$ and the asymptotic density from \eqref{eq: Asymptotic density from saddle points} for specific choices of $S$.

In the following, we derive the microscopic scaling limit at the origin.
We assume that the variance profile $S$ is irreducible, since otherwise the spectral measure of $H$ is a convex combination of the spectral measures of its irreducible components. Furthermore, from \cite{krueger2025SingularityDegree} we know that the asymptotic density of states has an isolated atom at zero whenever the variance profile $S$ does not have support. In this case the asymptotic eigenvalue statistics of $H$ in a vicinity of the origin is trivial. Therefore we assume in the following that $S$ has support, i.e.\ that there is a permutation matrix $P \in \{0, 1\}^{K \times K}$ such that $S \geq \eps P$ entry wise for some $\eps > 0$. In this case $\rho_{\infty}$ has a density on the real line that we will denote by the same symbol, i.e.\ $\rho_{\infty}(\dd E) = \rho_{\infty}(E) \dd E$.

From the classification of singularities of $\rho_{\infty}$ in \cite{krueger2025SingularityDegree} we know that the asymptotic density of states has a polynomial singularity at the origin such that the limit
\begin{equation} \label{eq: density limit constant}
    \vartheta := \lim_{E \to 0} \abs{E}^{\frac{\ell-1}{\ell+1}} \rho_{\infty}(E) \in (0, \infty)
\end{equation}
is finite and strictly positive, where $\ell \leq K$ is a positive integer that depends on the location of zero and non-zero entries within the variance profile $S$.
From the above relation \eqref{eq: density limit constant} we deduce that
\begin{equation} \label{eq: spacing scale}
    \lim_{N \to \infty} K N \int_{-\eta_N / 2}^{\eta_N / 2} \rho_{\infty}(E) \, \dd E = 1\,, \qquad
    \text{for } \eta_N := 2 \p1{\vartheta K N (\ell+1)}^{-\frac{\ell+1}{2}}.
\end{equation}
We interpret $\eta_N$ as the local spacing scale, i.e.\ the typical distance between consecutive eigenvalues at the origin.

For $K=1$ we have $\ell=1$ and $H$ is a GUE matrix, up to a constant.
In this case the spacing scale is proportional to $\frac{1}{N}$ and the local eigenvalue one-point function $\lim_{N \to \infty} K N \eta_N \rho_N(\eta_N x)$ at the origin is constant, i.e.\ independent of $x \in \R$.
For $S = 1$ we have
\begin{equation} \label{eq: K=1 expectation limit}
    \lim_{N \to \infty} \frac{1}{N} \E\s*{\Tr\p*{\frac{1}{H - N^{-1} \zeta}}} = \ii.
\end{equation}

For $K=2$ three cases occur.
When all entries of $S$ are non-zero the density $\rho_{\infty}$ remains bounded and $\ell=1$.
In this case the local eigenvalue statistics at the origin is given by the sine kernel process and the local eigenvalue one-point function is again constant.
Bounded density with $\ell=1$ also occurs when the diagonal entries of $S$ vanish while the off-diagonal entries are positive.
The corresponding matrix ensemble is called the \emph{chiral GUE} and was introduced in \cite{shuryak1993SumRules,verbaarschot1993SpectralDensity,verbaarschot1994SpectrumQCD} as a model for the Dirac operator in quantum chromodynamics.
Its eigenvalue distribution is a determinantal point process whose kernel can be written in terms of Laguerre polynomials.
The macroscopic limiting density is the semi-circle law, the microscopic scaling limit in the bulk is the sine kernel except at the origin where the limiting kernel can be written in terms of Bessel-$J$ functions \cite{forrester1993SpectrumEdge}.
Here the scaling limit of the one-point function is
\begin{equation} \label{eq: chiralGUE microscopic density}
    \lim_{N \to \infty} \rho_N\p*{\frac{\xi}{2 N}}
    = \frac{\abs{\xi}}{2} (J_0(\xi)^2 + J_1(\xi)^2).
\end{equation}
The remaining case with $K=2$ block indices is $\ell=2$, where up to permutation of the block indices the variance profile takes the form
\begin{equation} \label{eq: K=2 profile}
    S = \begin{pmatrix}
        s_{11} & s_{12} \\
        s_{12} & 0
    \end{pmatrix}, \qquad
    \text{with }\; s_{11}, s_{12} > 0.
\end{equation}
As we will see in Theorem~\ref{thm: microscopic resolvent expectation} below the eigenvalue one-point function on the local spacing scale $N^{-3/2}$ around the origin is non-trivial in this situation.

For $K=3$ block indices with a variance profile $S$ that has support there are only two cases. Either the variance profile is fully indecomposable (see \cite{krueger2025SingularityDegree} for the precise definition) or it is, up to permutation of the block indices, of the form
\begin{equation}\label{eq: K=3 profile}
    S = \begin{pmatrix}
        s_{11} & s_{12} & s_{13} \\
        s_{12} & s_{22} & 0 \\
        s_{13} & 0 & 0
    \end{pmatrix}, \qquad
    \text{with }\; s_{11} \geq 0 \;\text{ and }\; s_{12}, s_{13}, s_{22} > 0.
\end{equation}
In the former case, we have $\ell=1$ and the local eigenvalue one-point function is constant.
In the latter, $\ell=3$ with a local spacing scale of $N^{-2}$.

Our main theorem provides the universal scaling limits of the expectation of the resolvent for the two non-trivial cases \eqref{eq: K=2 profile} and \eqref{eq: K=3 profile}, the proof is given in Section~\ref{sec: Large N Limit at the Origin}.

\begin{theorem} \label{thm: microscopic resolvent expectation}
    Let $H$ be a Hermitian $K \times K$-block random matrix with centred independent Gaussian entries as in \eqref{eq: matrix model definition}.
    Let $\zeta \in \C \setminus \R$ with $\Im(\zeta) > 0$.

    For $K = 2$ and
    \[
        S = \begin{pmatrix}
            s_{11} & s_{12} \\
            s_{12} & 0
        \end{pmatrix}, \qquad
        \text{with }\; s_{11}, s_{12} > 0,
    \]
    we have
    \begin{equation} \label{eq: K=2 expectation limit}
        \begin{split}
        \lim_{N \to \infty} \frac{s_{12}}{\sqrt{s_{11}}} \frac{1}{N^{3/2}} \E\s*{\Tr\p*{\frac{1}{H - \frac{s_{12}}{\sqrt{s_{11}}} N^{-3/2} \zeta}}}
        & = \frac{\ii \sqrt{2}}{\sqrt{\pi}} \bigg[
            \frac{1}{2}
            \MeijerG{3, 0}{0, 3}{-}{0, \frac{1}{2}, \frac{1}{2}}{-\frac{\zeta^2}{8}}
            \HypergeometricPFQ{0}{2}{-}{\frac{1}{2}, 1}{-\frac{\zeta^2}{8}}
        \\
        & \hphantom{= \frac{\ii \sqrt{2}}{\sqrt{\pi}} \bigg[}
            + \MeijerG{3, 0}{0, 3}{-}{\frac{1}{2}, 1, \frac{3}{2}}{-\frac{\zeta^2}{8}}
            \HypergeometricPFQ{0}{2}{-}{\frac{3}{2}, 2}{-\frac{\zeta^2}{8}}
        \\
        & \hphantom{= \frac{\ii \sqrt{2}}{\sqrt{\pi}} \bigg[}
            + \MeijerG{3, 0}{0, 3}{-}{\frac{1}{2}, \frac{1}{2}, 1}{-\frac{\zeta^2}{8}}
            \HypergeometricPFQ{0}{2}{-}{1, \frac{3}{2}}{-\frac{\zeta^2}{8}}
        \bigg],
        \end{split}
    \end{equation}
    where $G_{p,q}^{m,n}$ is the Meijer $G$-function \cite[\DLMFsecref{16.17}]{NIST:DLMF} and ${}_p F_q$ is the generalized hypergeometric function \cite[\DLMFsecref{16.2}]{NIST:DLMF}.

    For $K = 3$ and
    \[
        S = \begin{pmatrix}
            s_{11} & s_{12} & s_{13} \\
            s_{12} & s_{22} & 0 \\
            s_{13} & 0 & 0
        \end{pmatrix}, \quad
        \text{with }\; s_{11} \geq 0 \;\text{ and }\; s_{12}, s_{13}, s_{22} > 0,
    \]
    we have
    \begin{eqsplit} \label{eq: K=3 expectation limit}
        \lim_{N \to \infty} \frac{s_{13} \sqrt{s_{22}}}{s_{12}} \frac{1}{N^2} \E\s*{\Tr\p*{\frac{1}{H - \frac{s_{13} \sqrt{s_{22}}}{s_{12}} N^{-2} \zeta}}}
        & = 2 \ii \bigg[
            K_0\p*{2 \sqrt{-\ii \zeta}} I_0\p*{2 \sqrt{-\ii \zeta}}
        \\
        & \hphantom{= 2 \ii \bigg[}
            + K_1\p*{2 \sqrt{-\ii \zeta}} I_1\p*{2 \sqrt{-\ii \zeta}}
        \bigg],
    \end{eqsplit}
    where $K_{\nu}$ and $I_{\nu}$ are the modified Bessel functions \cite[\DLMFsecref{10.25}]{NIST:DLMF}.
\end{theorem}

\begin{remark}
    The Meijer $G$-function also appears in the local one-point function for products of random matrices.
    As an example, consider the Hermitized product $G_m^* \dots G_1^* H G_1 \dots G_m$ of $m \in \N$ independent rectangular matrices $G_i$ with i.i.d.\ complex Gaussian entries and an independent GUE matrix $H$.
    It is known \cite{forrester2018MatrixProductEnsembles} that the eigenvalues form a determinantal point process, that the limiting density $\rho_{\infty}$ has a singularity at the origin and that, in the limit of large matrix size, the microscopic correlation kernel at the origin can be written in terms of the \emph{Meijer $G$-Kernel}
    \[
        K_{\text{Meijer}}^M(x, y) = \int_0^1 \dd u \, \MeijerG{1, 0}{0, M+1}{-}{-\nu_0, \dots, -\nu_M}{x u} \MeijerG{M, 0}{0, M+1}{-}{\nu_M, \dots, \nu_0}{y u},
    \]
    where $M = 2 m + 1$ is the number of matrices in the product and $\nu_{2 i}$ and $\nu_{2 i + 1}$ are related to the dimension of $G_i$.
    Furthermore, \cite[Prop.~13]{forrester2018MatrixProductEnsembles} gives a contour representation involving $G^{1, 0}_{0, m+1}(\dots) \times G^{m+1, 0}_{0, m+1}(\dots)$.
    Note that for $m = 2$ we can write $G^{1, 0}_{0, 3}$ as a ${}_{0} F_{2}$-function, see \cite[Sec.~6.4]{luke1969SpecialFunctionsApproximations}.
    But apart from these superficial similarities, we do not see a deeper connection between our $K = 2$ model \eqref{eq: matrix model definition} and the Hermitized product for $m = 2$ or for any other $m \in \N$.

    However, the $K = 2$ case may be related to a different product ensemble.
    Computing the characteristic polynomial via the Schur complement formula gives
    \begin{equation}
        \det\begin{pmatrix}
            A - \lambda & B \\
            B^* & 0 - \lambda
        \end{pmatrix}
        = \det\p*{(0 - \lambda) (A - \lambda - B (0 - \lambda)^{-1} B^*)}
        = \det(\lambda^2 - \lambda A + B B^*).
    \end{equation}
    For $\lambda$ close to zero we may naively assume that $\lambda^2$ is much smaller than $-\lambda A + B B^*$ and write
    \begin{equation}
        \det(\lambda^2 - \lambda A + B B^*)
        \approx \det(-\lambda A + B B^*)
        = \det(A) \det(- \lambda + A^{-1} B B^*)
        = \det(A) \det(B^* A^{-1} B - \lambda).
    \end{equation}
    Based on this naive calculation we conjecture that the origin limit for the $K = 2$ case agrees with the origin limit for the product $B^* A^{-1} B$, where $A$ is a GUE matrix and $B$ an independent complex Ginibre matrix.
\end{remark}

\begin{figure}[ht]
    \centering
    \begin{subfigure}[t]{0.4\textwidth}
        \centering
        \includegraphics[height=130pt]{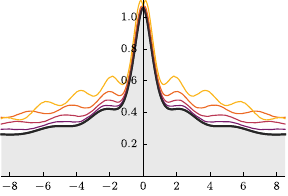}
        \caption{$K = 2$}
    \end{subfigure}
    \hfill
    \begin{subfigure}[t]{0.59\textwidth}
        \centering
        \includegraphics[height=130pt]{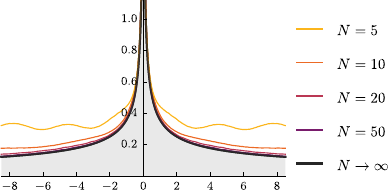}
        \caption{$K = 3$}
    \end{subfigure}
    \caption{Plots of $\xi \mapsto K N \eta_N \rho_N(\eta_N \xi)$, the microscopic scaling of the one-point function at the origin. The graphs for finite $N$ (coloured lines) are histograms obtained from numerical simulations of \eqref{eq: matrix model definition}. The limit $\lim_{N \to \infty} K N \eta_N \rho_N(\eta_N \xi)$ (black line) is derived from \eqref{eq: Stieltjes inversion} and Theorem~\ref{thm: microscopic resolvent expectation}.}
    \label{fig: microscopic one-point function}
\end{figure}

By using \eqref{eq: Stieltjes inversion} and \eqref{eq: K=2 expectation limit} (for $K = 2$) or \eqref{eq: K=3 expectation limit} (for $K = 3$) we obtain a closed formula for the local scaling limit at the origin, i.e.\ for $\xi \mapsto \lim_{N \to \infty} K N \eta_N \rho_N(\eta_N \xi)$.
We plot this function in Figure~\ref{fig: microscopic one-point function}.
More details are given in Section~\ref{sec: Microscopic Scaling Limit}.
In particular, we observe the following asymptotic behavior.

\begin{corollary} \label{thm: microscopic one-point function}
    As $\abs{\xi} \to 0$ we have
    \begin{equation} \label{eq: microscopic density origin}
        \lim_{N \to \infty} K N \eta_N \rho_N(\eta_N \xi) = \begin{cases}
            1, & K = 1, \\
            \frac{4 \pi}{3^{9/4}} + O(\abs{\xi}), & K = 2, \\
            -\frac{\pi}{4} \log(\abs{\xi}) + \frac{\pi}{2} \p*{\log\p*{\frac{2}{\pi}} - \gamma + \frac{1}{2}} + O(\abs{\xi}), & K = 3,
        \end{cases}
    \end{equation}
    where $\gamma \approx 0.5772$ is the Euler–Mascheroni constant and the local spacing scale $\eta_N$ is
    \begin{align}
        K & = 1: \quad \eta_N = \pi N^{-1}, \\
        K & = 2: \quad \eta_N = \frac{2 s_{12}}{\sqrt{s_{11}}} \p*{\frac{2 \pi}{3\sqrt{3}}}^{3/2} N^{-3/2}, \\
        K & = 3: \quad \eta_N = \frac{s_{13} \sqrt{s_{22}}}{s_{12}} \frac{\pi^2}{4} N^{-2}.
    \end{align}
    As $\abs{\xi} \to \infty$ we have the asymptotic
    \begin{equation} \label{eq: microscopic density infinity}
        \lim_{N \to \infty} K N \eta_N \rho_N(\eta_N \xi) = \begin{cases}
            1, & K = 1, \\
            \frac{2^{2/3}}{3} \abs{\xi}^{-1/3} + O(\abs{\xi}^{-1}), & K = 2, \\
            \frac{1}{2 \sqrt{2}} \abs{\xi}^{-1/2} + O(\abs{\xi}^{-1}), & K = 3.
        \end{cases}
    \end{equation}
\end{corollary}

\begin{remark}
    In the $\abs{\xi} \to 0$ asymptotics \eqref{eq: microscopic density origin} for $K = 3$, we see that the microscopic one-point function has a logarithmic singularity and we may interpret this singularity as an attraction of the eigenvalues to the origin.
    For comparison, let us consider a 1-dimensional log-gas in an external potential $V$ with a local attraction due to the insertion of a point charge (also called a root-type Fisher-Hartwig singularity).
    The corresponding unitary ensemble of random matrices has a probability density proportional to $\abs{\det(M)}^{2 \alpha} \ee^{-N \Tr(V(M))} \dd M$, where $M$ is a Hermitian $N \times N$-matrix and $\alpha > -\frac{1}{2}$ is the strength of the charge inserted at the origin.
    If $\alpha < 0$, the charge is attractive.
    Note that the point charge affects the eigenvalue distribution only on the microscopic scale, because $\alpha$ does not scale with $N$.
    It is known from \cite{kuijlaars2003UniversalityOriginSpectrum} that the microscopic scaling limit at the origin is universal for a large class of potentials $V$.
    However, one of the assumptions on $V$ is that the equilibrium measure $\rho_{\infty}$ is positive and finite at the origin, unlike our model.
    The limiting correlation kernel is given in terms of Bessel-$J$ functions (but is different from the hard edge Bessel-kernel) and for $\abs{\xi} \to 0$ the microscopic one-point function is asymptotically proportional to $\abs{\xi}^{\alpha}$.
    Hence, there is a divergence but not a logarithmic divergence when $\alpha < 0$.
    This comparison suggests that the local attraction in our case is weaker than the attraction due to a charge of opposite sign.

    The asymptotic behaviour for $\abs{\xi} \to \infty$ in \eqref{eq: microscopic density infinity} can be explained easily: we transition from the microscopic scale to the mesoscopic scale, on which the macroscopic singularity $\abs{\xi}^{-\frac{l-1}{l+1}}$ of $\rho_{\infty}$ dictates the asymptotics.
\end{remark}

\begin{remark}
    We observe that the spacing scale $\eta_N$ in the $K = 3$ case does not depend on the value of $s_{11}$.
    This is consistent with the idea from \cite{krueger2025SingularityDegree} that the singularity degree of $\rho_{\infty}$ only depends on the length of a certain path in the 0-1 mask of $S$ and that the limit $\vartheta$ from \eqref{eq: density limit constant} depends only on the entries of $S$ on this path.
    For the variance profile in \eqref{eq: K=3 profile} the path does not include the $s_{11}$-entry and therefore it does not influence the singularity degree or the limit $\vartheta$ (compare \eqref{eq: K=3 spacing scale}) which we use to define the spacing scale in \eqref{eq: spacing scale}.
\end{remark}

Comparing the microscopic scaling limit \eqref{eq: chiralGUE microscopic density} for the chiral GUE with our $K = 2$ result, we see that for $s_{11} = 0$ the local one-point function contains Bessel-$J$ functions while for $s_{11} > 0$ it involves Meijer $G$ and generalized hypergeometric functions.
To analyze the transition between the two formulas, we consider the large-$N$ limit where $s_{11} \to 0$ as $N \to \infty$.
We call this regime the \emph{weak non-chirality limit}.
For the analogous limit in the $K = 3$ case, we take $s_{11} = 0$ and $s_{12} = s_{21} \to 0$ as $N \to \infty$.

We mention that a different weak non-chirality limit of the chiral GUE, namely when both diagonal elements of $S$ are non-zero and $s_{11} = s_{22} \to 0$ as $N \to \infty$, was the topic of previous work \cite{damgaard2010MicroscopicSpectrum, akemann2011WilsonDirac}.

\begin{figure}[ht]
    \centering
    \begin{subfigure}[t]{0.49\textwidth}
        \centering
        \includegraphics[height=160pt]{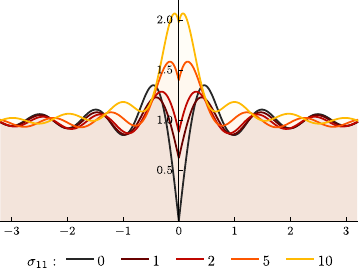}
        \caption{$K = 2$}
    \end{subfigure}
    \hfill
    \begin{subfigure}[t]{0.49\textwidth}
        \centering
        \includegraphics[height=160pt]{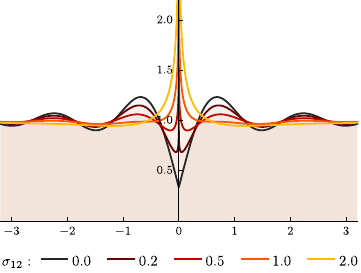}
        \caption{$K = 3$}
    \end{subfigure}
    \caption{Plots of $\xi \mapsto \lim_{N \to \infty} K N \eta_N \rho_N(\eta_N \xi)$ in the weak non-chirality limit derived from Theorem~\ref{thm: microscopic resolvent expectation weak non-chirality}.}
    \label{fig: microscopic one-point function weak non-chirality}
\end{figure}

\begin{theorem} \label{thm: microscopic resolvent expectation weak non-chirality}
    Let $H$ be a $2N \times 2N$ block matrix as in \eqref{eq: matrix model definition} with
    \[
        S = \begin{pmatrix}
            N^{-1} \sigma_{11} & 1 \\
            1 & 0 \\
        \end{pmatrix}, \qquad
        \text{for }\; \sigma_{11} \in [0, \infty).
    \]
    Let $\zeta \in \C$ with $\Im(\zeta) > 0$, then
    \begin{eqsplit}
        \lim_{N \to \infty} \frac{1}{N} \E\s*{\Tr\p*{\frac{1}{H - N^{-1} \zeta}}}
        & = \frac{\ii}{2} \sum_{n, r \in \Z} c_{n, r} \int_0^{\infty} \dd v \, \ee^{-\p*{\frac{1}{2} \sigma_{11} v^2 - \ii \zeta (v + \frac{1}{v})}} v^n
        \\
        & \qquad
            \times \sum_{k, l = 0}^{\infty} \frac{(\sigma_{11} / 2)^k (-\ii \zeta)^{2k + 2l + 1 + r}}{k! \, l! \, (2 k + l + 1 + r)!},
    \end{eqsplit}
    where the only non-zero coefficients $c_{n, r}$ are
    \begin{align*}
        & c_{-3, -1} = c_{-2, -2} = c_{-2, 0} = c_{0, -2} = c_{0, 0} = -\ii \zeta, \qquad
        c_{-1, -1} = -1 - 2 \ii \zeta, \qquad
        c_{1, -1} = -1 - \ii \zeta, \\
        & c_{-2, 1} = c_{0, -1} = c_{0, 1} = c_{2, -1} = \sigma_{11}, \qquad
        c_{-1, 0} = c_{1, 0} = 2 \sigma_{11}.
    \end{align*}

    For the $K = 3$ case, let $H$ be a $3N \times 3N$ block matrix as in \eqref{eq: matrix model definition} with
    \[
        S = \begin{pmatrix}
            0 & N^{-1} \sigma_{12} & 1 \\
            N^{-1} \sigma_{12} & 1 & 0 \\
            1 & 0 & 0
        \end{pmatrix}, \qquad
        \text{for }\; \sigma_{12} \in [0, \infty).
    \]
    Let $\zeta \in \C$ with $\Im(\zeta) > 0$, then
    \begin{equation}
    \begin{split}
        \lim_{N \to \infty} \frac{1}{N} \E\s*{\Tr\p*{\frac{1}{H - N^{-1} \zeta}}}
        & = 2 \ii \Big[ (\sigma_{12} - 2 \ii \zeta) \p1{K_0(x) I_0(x) + K_1(x) I_1(x)}
        \\
        & \hphantom{= 2 \ii \Big[}
            + \sqrt{-\ii \zeta (\sigma_{12} - \ii \zeta)} \p1{K_1(x) I_0(x) + K_0(x) I_1(x)} \Big],
    \end{split}
    \end{equation}
    where $x := 2 \sqrt{-\ii \zeta (\sigma_{12} - \ii \zeta)}$.
\end{theorem}
Note that for $K = 2$ the condition $\Im(\zeta) > 0$ is necessary for the convergence of the integrals.
The formula for $K = 3$ also holds for $\zeta \in \R$, but when $\sigma_{12} = 0$ the choice of the square root branch must be consistent with the $\Im(\zeta) \searrow 0$ limit.
For $K = 2$ with $\sigma_{11} = 0$ and $K = 3$ with $\sigma_{12} = 0$ the macroscopic limiting density is the semi-circle law with radius $2$.
Hence the local spacing scale at the origin is $\eta_N = \frac{\pi}{K N}$.
In Figure~\ref{fig: microscopic one-point function weak non-chirality} we plot the resulting microscopic scaling limit of the one-point functions.

\begin{remark}
    In the case $K = 2$ with $\sigma_{11} = 0$, the $v$-integrals evaluate to Bessel-$K$ and the double-sums to Bessel-$I$ functions, such that
    \begin{equation} \label{eq: chiralGUE expectation of resolvent}
        \lim_{N \to \infty} \frac{1}{N} \E\s*{\Tr\p*{\frac{1}{H - N^{-1} \zeta}}}
        = 4 \zeta \p*{K_0(-2 \ii \zeta) I_0(-2 \ii \zeta) + K_1(-2 \ii \zeta) I_1(-2 \ii \zeta)}.
    \end{equation}
    We recover the local one-point function in \eqref{eq: chiralGUE microscopic density} by using the relations \cite[\DLMFeqref{10.27}{6} and \DLMFeqref{10.27}{11}]{NIST:DLMF} and the fact that $\lim_{\eps \searrow 0} \Im(Y_n(x + \ii \eps)) = 2 J_n(x) \theta(-x)$ for $n \in \Z$ and $x \in \R$, where $\theta$ is the Heaviside function.

    In the case $K = 3$, a series expansion of the Bessel functions shows that for $\abs{\xi} \to 0$ we have
    \[
        \lim_{N \to \infty} K N \eta_N \rho_N\p*{\eta_N \xi}
        \sim - \frac{\sigma_{12}}{3} \p*{ \log(\abs{\xi}) + \log\p*{\frac{\pi}{3} \sigma_{12}} + 2 \gamma - 1 } + \frac{1}{3},
    \]
    where $\gamma$ is again the Euler-Mascheroni constant.
    The logarithmic divergence noted in Corollary~\ref{thm: microscopic one-point function} appears as soon as $\sigma_{12} > 0$.
    For $\sigma_{12} = 0$ we get the local one-point function
    \[
        \lim_{N \to \infty} \rho_N\p*{\frac{\xi}{3 N}}
        = \frac{2}{3} \abs{\xi} \p*{J_0\p*{\frac{2}{3} \xi}^2 + J_1\p*{\frac{2}{3} \xi}^2} + \frac{1}{3 \pi}.
    \]
    This combination of $\frac{2}{3}$ times the chiral GUE one-point function plus $\frac{1}{3}$ times a constant is expected, since an eigenvalue of the matrix
    \[
        H = \begin{pmatrix}
            0 & 0 & A \\
            0 & B & 0 \\
            A^* & 0 & 0
        \end{pmatrix}
    \]
    is either an eigenvalue of $\begin{pmatrix} 0 & A \\ A^* & 0 \end{pmatrix}$ or an eigenvalue of $B$.
    These blocks also correspond to the irreducible components of $S$.
    In our model, \eqref{eq: matrix model definition} $B$ is a GUE matrix and so its local one-point function at the origin is constant.
\end{remark}

\section{Integral Representation for Finite \texorpdfstring{$N$}{N}}
\label{sec: Integral Representation}
In this section we prove Lemma~\ref{thm: Finite N Integral Representation}.
We use the supersymmetric method and the superbosonization formula similar to the example in \cite{bunder2007SuperbosonizationFormula}, but our notation is closer to \cite{cipolloni2020OptimalLower} so that we can apply the superbosonization formula from \cite{littelmann2008SuperbosonizationInvariant}.

\subsection{SUSY Representation}
We start with the identity
\begin{equation} \label{eq: matrix inverse as commuting integral}
    A^{-1} = \det(A) \int_{\C^n} [Dx] \, \ee^{-x^* A x} x x^*\,, \qquad
    [Dx] := \prod_i \frac{\dd\!\Re x_i \,\dd\!\Im x_i}{\pi}
\end{equation}
which holds for any $n \in \N$ and $A \in \C^{n \times n}$ such that $\Re \spec(A) \subset (0, \infty)$, where $\spec(A)$ is the set of eigenvalues of $A$.
Note that we write $x^*$ for the row vector of complex conjugates $\ol{x_1}, \dots, \ol{x_n}$ and that we multiply vectors and matrices row-by-column.
In particular, $x^* A x = \sum_{i, j = 1}^n \ol{x_i} A_{ij} x_j$ is a $1 \times 1$-matrix which we interpret as a scalar and $x x^*$ is the $n \times n$-matrix with entries $(x x^*)_{ij} = x_i \ol{x_j}$ for $i, j = 1, \dots, n$.

Next, we express the determinant of $A$ using \emph{Grassmann variables}.
Let $\xi_i$ and $\ol{\xi}_i$ for $i = 1, \dots, n$ be anti-commuting variables that satisfy
\[
    \xi_j \xi_i = -\xi_i \xi_j, \qquad
    \ol{\xi}_j \xi_i = -\xi_i \ol{\xi}_j, \qquad
    \ol{\xi}_j \ol{\xi}_i = -\ol{\xi}_i \ol{\xi}_j \qquad
    \text{for all } i, j = 1, \dots, n.
\]
In particular, $\xi_i^2 = \ol{\xi}_i^2 = 0$.
Additionally, we introduce the derivatives $\partial_{\xi_i}$ and $\partial_{\ol{\xi}_i}$ as
\[
    \partial_{\xi_i} \xi_j = \partial_{\ol{\xi}_i} \ol{\xi}_j = \delta_{i, j}, \qquad
    \partial_{\xi_i} \ol{\xi}_j = \partial_{\ol{\xi}_i} \xi_j = 0, \qquad
    \partial_{\xi_i} 1 = \partial_{\ol{\xi}_i} 1 = 0.
\]
The determinant identity is then (see e.g.\ \cite[Sec.~5.13]{haake2018QuantumSignatures}, but note that we use a different sign convention)
\begin{equation} \label{eq: Grassmann for determinant}
    \det(A) = [\partial_{\xi} \partial_{\ol{\xi}}] \, \ee^{\xi^* A \xi}\,, \qquad
    [\partial_{\xi} \partial_{\ol{\xi}}] := \partial_{\xi_1} \partial_{\ol{\xi}_1} \dots \partial_{\xi_n} \partial_{\ol{\xi}_n}\,,
\end{equation}
where the exponential of anti-commuting variables is defined by its (terminating) Taylor series, $\xi = (\xi_1, \dots, \xi_n)^{\top}$ is a column vector while $\xi^* = (\ol{\xi}_1, \dots, \ol{\xi}_n)$ is a row vector and $\xi^* A \xi = \sum_{i, j = 1}^n \ol{\xi}_i A_{ij} \xi_j$.

When combining the identities \eqref{eq: matrix inverse as commuting integral} and \eqref{eq: Grassmann for determinant}, it is useful to introduce the $n \times (1|1)$-matrix $\Psi = (x, \xi)$ and the $(1|1) \times n$-matrix $\Psi^* = \begin{pmatrix}x^* \\ \xi^*\end{pmatrix}$.
The dimension $(p|q)$ denotes the presence of $p \in \N$ commuting and $q \in \N$ anti-commuting components.
Then
\begin{equation}
    x^* A x - \xi^* A \xi
    = \Str\begin{pmatrix}
        x^* A x & x^* A \xi \\
        \xi^* A x & \xi^* A \xi
    \end{pmatrix}
    = \Str(\Psi^* A \Psi),
\end{equation}
where the \emph{supertrace} $\Str$ of a $(p|q) \times (p|q)$-supermatrix is defined as the scalar
\[
    \Str\begin{pmatrix}
        a & \sigma \\
        \tau & b
    \end{pmatrix}
    := \Tr(a) - \Tr(b).
\]
The $p \times p$-block $a$ is usually called the \emph{boson-boson} block and the $q \times q$-block $b$ is called the \emph{fermion-fermion} block, both have commuting entries.
The $p \times q$-block $\sigma$ and the $q \times p$-block $\tau$ contain Grassmann variables.
Observe that we only need the case $p = q = 1$.

Using the cyclic invariance of the trace and the anti-commutativity of $\xi$ and $\xi^*$ yields
\begin{equation}
    x^* A x - \xi^* A \xi
    = \Tr(A x x^*) + \Tr(A \xi \xi^*)
    = \Tr(A \Psi \Psi^*).
\end{equation}
We conclude that $\Str(\Psi^* A \Psi) = \Tr(A \Psi \Psi^*)$, where $\Psi^* A \Psi$ is a $(1|1) \times (1|1)$-supermatrix while $\Psi \Psi^* = x x^* + \xi \xi^*$ is a $n \times n$-matrix whose entries contain linear combinations of commuting variables that are either complex numbers or quadratic expressions of Grassmann variables.

To summarize, for $A \in \C^{n \times n}$ with $\Re \spec(A) \subset (0, \infty)$ we obtain the formula
\begin{equation}
    \Tr\p*{A^{-1}}
    = \int_{\C^n} [Dx] [\partial_{\xi} \partial_{\ol{\xi}}] \, \ee^{-x^* A x + \xi^* A \xi} \Tr(x x^*)
    = \int_{\C^n} [Dx] [\partial_{\xi} \partial_{\ol{\xi}}] \, \ee^{-\Str(\Psi^* A \Psi)} x^* x
\end{equation}
and in particular, for $A = \ii (H - z)$ with $\Im(z) > 0$ and $n = K N$, we get
\begin{equation}\label{eq: single resolvent as integral}
    \Tr\p*{\frac{1}{H-z}}
    = \ii \Tr(A^{-1})
    = \ii \int_{\C^{KN}} [Dx] [\partial_{\xi} \partial_{\ol{\xi}}] \, \ee^{-\ii \Tr(H \Psi \Psi^*) + \ii z \Str(\Psi^* \Psi)} x^* x \,.
\end{equation}
Note that when defining $A$ we include the factor $\ii$ to ensure that the condition $\Re \spec(A) = \Im(z) > 0$ holds.

\subsection{Evaluation of the Expectation}
The expectation value
\[
    \E\s*{\ee^{-\ii \Tr(H \Psi \Psi^*)}}
    = \E\s*{\ee^{-\ii \Tr(H Y)}}, \qquad
    Y := \Psi \Psi^* = x x^* + \xi \xi^*
\]
is the Fourier-transform of the matrix probability distribution.
Since $H$ contains only centred Gaussian random variables, the expectation satisfies
\begin{equation}
    \E\s*{\ee^{-\ii \Tr(H Y)}}
    = \exp\p*{-\frac{1}{2} \E\s*{(\Tr H Y)^2}}.
\end{equation}
Since $H$ is Hermitian, the entries $h_{\alpha \beta}$ and $h_{\gamma \delta}$ are independent except for the index pairs $(\alpha, \beta) = (\gamma, \delta)$ and $(\alpha, \beta) = (\delta, \gamma)$.
Hence
\begin{equation}
    \E\s*{(\Tr H Y)^2}
    = \sum_{\alpha, \beta, \gamma, \delta = 1}^{K N} \E\s*{h_{\alpha \beta} Y_{\beta \alpha} h_{\gamma \delta} Y_{\delta \gamma}}
    = \sum_{\alpha, \beta, \gamma, \delta = 1}^{K N} \p*{\delta_{\alpha \gamma} \delta_{\beta \delta} \E\s*{h_{\alpha \beta}^2} + \delta_{\alpha \delta} \delta_{\beta \gamma} \E\s*{\abs{h_{\alpha \beta}}^2}} Y_{\beta \alpha} Y_{\delta \gamma}
\end{equation}
The blocks $X_{ij}$ in our matrix model \eqref{eq: matrix model definition} are drawn from the complex Ginibre ensemble and therefore we have $\E\1 h_{\alpha \beta}^2 = 0$.
The remaining expectation $\E\1 \abs{h_{\alpha \beta}}^2$ is equal to $\frac{1}{N} s_{ij}$ when $\alpha$ belongs to the $i$-th and $\beta$ to the $j$-th batch of indices, i.e.\ when $N (i - 1) < \alpha \leq N i$ and $N (j - 1) < \beta \leq N j$.
To keep the formulas short, we introduce $E_N$ as the $N \times N$-matrix that contains $1$ in every entry and write $\E\1 \abs{h_{\alpha \beta}}^2 = (\frac{1}{N} S \otimes E_N)_{\alpha \beta}$.
Thus
\[
    \E\s*{\ee^{-\ii \Tr(H \Psi \Psi^*)}}
    = \exp\p*{-\frac{1}{2 N} \sum_{\alpha, \beta = 1}^{K N} (S \otimes E_N)_{\alpha \beta} Y_{\beta \alpha} Y_{\alpha \beta}}, \qquad
    \text{where } Y_{\alpha \beta} = \Psi_{\alpha} \Psi_{\beta}^*
    = (x_{\alpha}, \xi_{\alpha}) \begin{pmatrix}
        \ol{x_{\beta}} \\
        \ol{\xi}_{\beta}
    \end{pmatrix}.
\]
Note that the supertrace is cyclic which allows us to write $Y_{\beta \alpha} Y_{\alpha \beta} = \Psi_{\beta} \Psi_{\alpha}^* \Psi_{\alpha} \Psi_{\beta}^* = \Str(\Psi_{\alpha}^* \Psi_{\alpha} \Psi_{\beta}^* \Psi_{\beta})$.
We now replace the Greek index $\alpha = 1, \dots, K N$ by a composite of two Latin indices $\alpha = (i, l)$, where $i = 1, \dots, K$ enumerates the blocks in our model \eqref{eq: matrix model definition} and $l = 1, \dots, N$ enumerates the entries inside the block:
\begin{equation}
    \E\s*{\ee^{-\ii \Tr(H \Psi \Psi^*)}}
    = \exp\p*{
        -\frac{1}{2 N} \sum_{i, j = 1}^{K} s_{ij} \Str\s*{
            \sum_{l = 1}^N \Psi_{(i, l)}^* \Psi_{(i, l)} \sum_{m = 1}^N \Psi_{(j, m)}^* \Psi_{(j, m)}
        }
    }.
\end{equation}
After we define the supermatrices $Q_i := \sum_{l = 1}^N \Psi_{(i, l)}^* \Psi_{(i, l)}$ for $i = 1, \dots, K$, we arrive at
\begin{equation}
    \E\s*{\Tr\p*{\frac{1}{H - z}}}
    = \ii \int_{\C^{KN}} [Dx] [\partial_{\xi} \partial_{\ol{\xi}}] \, \exp\p*{-\frac{1}{2 N} \sum_{i, j = 1}^{K} s_{ij} \Str(Q_i Q_j) + \ii z \sum_{i = 1}^K \Str(Q_i)} x^* x \,.
\end{equation}
Observe that $x^* x$ is the boson-boson block of $Q_i$ summed over $i = 1, \dots, K$.

\subsection{Superbosonization}
The matrix $Q_i$ only depends on inner products of $(x_{(i, l)})_{l = 1}^N$, $(\xi_{(i, l)})_{l = 1}^N$ and $(\ol{\xi}_{(i, l)})_{l = 1}^N$.
The remaining degrees of freedom can be integrated out using the \emph{superbosonization formula}. The version from \cite{littelmann2008SuperbosonizationInvariant} for $U(n)$-symmetry and $p = q = 1$ reads
\begin{equation} \label{eq: Superbosonization formula}
    \int_{\C^n} [D Z] [\partial_{\zeta} \partial_{\ol{\zeta}}] \, F\begin{pmatrix}
        \wt{Z} Z & \wt{Z} \zeta \\
        \wt{\zeta} Z & \wt{\zeta} \zeta
    \end{pmatrix}
    = \int_D [DQ] \, \Sdet^N(Q) F(Q), \qquad
    Q = \begin{pmatrix}
        a & \sigma \\
        \tau & b
    \end{pmatrix}.
\end{equation}
On the left-hand side we have $Z \in \C^n$, $\wt{Z} = (\ol{Z_1}, \dots, \ol{Z_n})$, $\zeta = (\zeta_1, \dots, \zeta_n)^\top$ and $\wt{\zeta} = (\ol{\zeta}_1, \dots, \ol{\zeta}_n)$, where
$\zeta_i$ and $\ol{\zeta}_i$ are independent anti-commuting variables.
On the right-hand side the integration domain is $D = (0, \infty) \times U(1)$, i.e.\ $a$ is a real positive number and $b$ is a complex number on the unit circle.
The symbols $\sigma$ and $\tau$ are anti-commuting variables and the integration measure is
\[
    [D Q] = \frac{\dd a}{a} \, \dd\mathrm{Vol}(b) \frac{1}{2 \pi} \partial_{\sigma} \partial_{\tau} \underbrace{(a - \sigma b^{-1} \tau) (b - \tau a^{-1} \sigma)}_{= a b - \tau \sigma - \sigma \tau + 0 = a b}
    = \frac{1}{2 \pi} \dd a \, \dd\mathrm{Vol}(b) \, b \, \partial_{\sigma} \partial_{\tau}.
\]
The measure $\dd\mathrm{Vol}(b)$ is the Haar measure on $U(1)$.
The integrand on the right-hand side of \eqref{eq: Superbosonization formula} also includes the \emph{superdeterminant} of $Q$, which in the case of a $(1|1) \times (1|1)$-supermatrix is simply given by
\[
    \Sdet(Q) = \frac{a}{b - \tau a^{-1} \sigma} = \frac{a}{b} \p*{1 - \frac{\tau \sigma}{a b}}^{-1}
\]
and the inverse on the right is defined by the Taylor series of $(1 - t)^{-1}$ around $t = 0$.

We use the superbosonization formula with $n = N$ and apply it $K$ times, with $Z = (x_{(i, l)})_{l = 1}^N$, $\zeta = (\xi_{(i, l)})_{l = 1}^N$ and $\wt{\zeta} = (\ol{\xi}_{(i, l)})_{l = 1}^N$ for $i = 1, \dots, K$, and arrive at
\begin{equation}
    \E\s*{\Tr\p*{\frac{1}{H - z}}}
    = \frac{\ii}{(2 \pi)^K} \int_{(0, \infty)^K} \dd a \int_{U(1)^K} \dd\mathrm{Vol}(b) \, \p*{\prod_{i = 1}^K b_i} [\partial_{\sigma} \partial_{\tau}] \, \ee^{f(Q)} \sum_{i = 1}^K a_i,
\end{equation}
where $\dd a = \prod_{i = 1}^K \dd a_i$, $\dd\mathrm{Vol}(b) = \prod_{i = 1}^K \dd\mathrm{Vol}(b_i)$, $[\partial_{\sigma} \partial_{\tau}] = \prod_{i = 1}^K \partial_{\sigma_i} \partial_{\tau_i}$ and the exponent is
\begin{equation}
    f(Q) = - \frac{1}{2 N} \sum_{i, j = 1}^K s_{ij} \Str(Q_i Q_j) + \ii z \sum_{i = 1}^K \Str(Q_i) + N \sum_{i = 1}^K \log \Sdet(Q_i), \qquad
    Q_i = \begin{pmatrix}
        a_i & \sigma_i \\
        \tau_i & b_i
    \end{pmatrix}.
\end{equation}

\subsection{Anti-Commuting Differentiation}
From the supertraces we get
\[
    \Str(Q_i) = a_i - b_i, \qquad
    \Str(Q_i Q_j) = a_i a_j + \sigma_i \tau_j - (\tau_i \sigma_j + b_i b_j)
    = a_i a_j - b_i b_j - \tau_j \sigma_i - \tau_i \sigma_j.
\]
For $\log \Sdet(Q_i)$ we use the Taylor expansion $\log(1 - x) = 1 - x + O(x^2)$ to compute
\[
    \log \Sdet(Q_i)
    = \log(a_i) - \log(b_i) - \log\p*{1 - \frac{\tau_i \sigma_i}{a_i b_i}}
    = \log(a_i) - \log(b_i) + \frac{\tau_i \sigma_i}{a_i b_i}.
\]
Hence the exponent is
\begin{eqsplit*}
    f(Q) & = - \frac{1}{2 N} \sum_{i,j} s_{ij} a_j a_i + \ii z \sum_{i} a_i + N \sum_{i} \log(a_i)
    + \frac{1}{2 N} \sum_{i,j} s_{ij} b_j b_i - \ii z \sum_{i} b_i - N \sum_{i} \log(b_i)
    \\
    & \quad + \frac{1}{2 N} \sum_{i,j} s_{ij} \tau_j \sigma_i + \frac{1}{2 N} \sum_{i,j} s_{ij} \tau_i \sigma_j + N \sum_{i} \frac{\tau_i \sigma_i}{a_i b_i}
    \\
    & = -N I\p*{\frac{a}{N}} + N I\p*{\frac{b}{N}} + \sum_{i,j} \tau_j \s2{\frac{s_{ji} + s_{ij}}{2 N} + \frac{N}{a_i b_i} \delta_{ij}} \sigma_i,
\end{eqsplit*}
where the function $I$ is defined in \eqref{eq: I definition}.
Note that $S$ is symmetric by construction and therefore $\frac{s_{ji} + s_{ij}}{2 N} = \frac{1}{N} s_{ij}$.
We now apply the identity \eqref{eq: Grassmann for determinant} again, with $\xi^*$ replaced by $\tau^\top = (\tau_1, \dots, \tau_K)$ and $\xi$ replaced by $\sigma = (\sigma_1, \dots, \sigma_K)^\top$, to obtain
\begin{eqsplit*}
    [\partial_{\sigma} \partial_{\tau}] \, \ee^{f(Q)}
    & = \ee^{-N I\p*{\frac{a}{N}} + N I\p*{\frac{b}{N}}} [\partial_{\sigma} \partial_{\tau}] \, \exp\p*{
        \tau^\top \s*{ \frac{1}{N} S + N \diag_{i = 1}^K\p*{\frac{1}{a_i b_i}} } \sigma
    }
    \\
    & = \ee^{-N I\p*{\frac{a}{N}} + N I\p*{\frac{b}{N}}} \det\s*{ \frac{1}{N} S + \diag_{i = 1}^K\p*{\frac{N}{a_i b_i}} }.
\end{eqsplit*}

To derive \eqref{eq: resolvent expectation integral}, we perform two more steps.
First, we write the integral over the Haar measure $\dd\mathrm{Vol}(b_i)$ as a counter-clockwise contour integral over the unit circle around $0$: $\int_{U(1)} \dd\mathrm{Vol}(b_i) \, b_i = \frac{1}{\ii} \oint \dd b_i$.
Second, we perform the coordinate transformation $a_i \to N a_i$ and $b_i \to N b_i$:
\begin{align*}
    \E\s*{\Tr\p*{\frac{1}{H - z}}}
    & = \frac{\ii}{(2 \pi \ii)^K} \int_{(0, \infty)^K} \dd a \, N^K \oint^K \, \dd b \, N^K \ee^{-N I(a) + N I(b)} \det\s*{ \frac{1}{N} S + \diag_{i = 1}^K\p*{\frac{N}{N^2 a_i b_i}} } \sum_{j = 1}^K N a_j
    \\
    & = \frac{\ii N^{2 K - K + 1}}{(2 \pi \ii)^K} \int_{(0, \infty)^K} \dd a \oint^K \, \dd b \, \ee^{-N I(a) + N I(b)} \det\s*{ S + \diag_{i = 1}^K\p*{\frac{1}{a_i b_i}} } \sum_{j = 1}^K a_j.
\end{align*}

\section{Large \texorpdfstring{$N$}{N} Limit for \texorpdfstring{$z \to 0$}{z -> 0}}
\label{sec: Large N Limit at the Origin}
In this section we prove Theorem~\ref{thm: microscopic resolvent expectation}.
As $N \to \infty$ the largest contributions to the integral in \eqref{eq: resolvent expectation integral} comes from regions around the critical points $\nabla I(a) = 0$ and $\nabla I(b) = 0$.
Therefore we use the saddle-point method (summarized in Appendix~\ref{app: Saddle-Point Method}) to approximate the integral.
For $z \neq 0$ fixed this approach leads to the known macroscopic density \eqref{eq: Asymptotic density from saddle points}.
We do not perform this calculation since the final result can be inferred from \cite{ajanki2017SingularitiesSolutions}.

To derive the microscopic scaling limit at the origin, we write $z = \eta \zeta$ with the local coordinate $\zeta \in \C$, $\Im(\zeta) > 0$ and the local spacing scale $\eta \equiv \eta_N \asymp N^{-\frac{\ell+1}{2}}$.
The notation $a_N \asymp b_N$ means that $a_N$ and $b_N$ are of the same order, i.e.\ there exist $N$-independent constants $0 < c \leq C < \infty$ such that $c \abs{b_N} \leq \abs{a_N} \leq C \abs{b_N}$ for sufficiently large $N$.

\subsection{\texorpdfstring{$2 \times 2$}{2x2} case}
Let $S$ be of the form
\[
    S = \begin{pmatrix}
        s_{11} & s_{12} \\
        s_{12} & 0
    \end{pmatrix}, \qquad
    \text{with }\; s_{11}, s_{12} > 0.
\]
The saddle-points of the $a$-integral in \eqref{eq: resolvent expectation integral} are the points $a \in \C^2$ such that $\nabla I(a) = 0$, i.e.\ the coordinates $a_1$ and $a_2$ satisfy
\begin{equation} \label{eq: K=2 saddle-point equation}
    1 = s_{11} a_1^2 + \frac{s_{12} a_1}{s_{12} a_1 - \ii z} - \ii z a_1, \qquad
    a_2 = \frac{1}{s_{12} a_1 - \ii z}.
\end{equation}
The same equations hold for the saddle-points of the $b$-integral, with $a_j$ replaced by $b_j$.
For $z = \eta \zeta$ with $\eta \to 0$ we see that $\abs{a_1} \asymp \eta^{1/3}$ and $\abs{a_2} \asymp \eta^{-1/3}$.
We refer to Section~\ref{sec: Weak Non-Chirality Limit for K = 2} for a procedure to determine these rates.

\subsubsection{Transformation of the integrals}
As $\eta \to 0$ the saddle-point has a degenerate limit: $a_1 \to 0$ and $\abs{a_2} \to \infty$.
To correct for this behaviour, we introduce the rescaled coordinates $\wt{a}_1 := \eta^{-1/3} a_1$ and $\wt{a}_2 := \eta^{1/3} a_2$ and similarly for the $b$ coordinates.
Thus, the $I(a)$-term becomes
\begin{eqsplit*}
    I(a)
    & = \frac{1}{2} \p*{\eta^{2/3} s_{11} \wt{a}_1^2 + \eta^{0} s_{12} \wt{a}_2 \wt{a}_1 + \eta^{0} s_{21} \wt{a}_1 \wt{a}_2}
    - \ii \eta \zeta \p*{\eta^{1/3} \wt{a}_1 + \eta^{-1/3} \wt{a}_2}
    - \p*{\log(\eta^{1/3} \wt{a}_1) + \log(\eta^{-1/3} \wt{a}_2)}
    \\
    & = s_{12} \wt{a}_1 \wt{a}_2 - \log(\wt{a}_1) - \log(\wt{a}_2) + O\p*{\eta^{2/3}} \qquad \text{for } \eta \to 0.
\end{eqsplit*}
In the $\eta \to 0$ limit, the saddle-points of the rescaled integral satisfy
\[
    0 = s_{12} \wt{a}_2 - \frac{1}{\wt{a}_1}
    \quad \text{and} \quad
    0 = s_{12} \wt{a}_1 - \frac{1}{\wt{a}_2}.
\]
However, this system of equations has an infinite number of solutions: $\wt{a}_2 = (s_{12} \wt{a}_1)^{-1}$.
To address this problem, we transform $(\wt{a}_1, \wt{a}_2)$ to new variables $(u_a, v)$ such that the saddle-line becomes a saddle-point for $u_a$ and the action is independent of $v$.
We write $\wt{a}_1 := u_a v$, $\wt{a}_2 := \frac{u_a}{v}$ for $u_a > 0$ and $v > 0$, then the $u_a$-integral has a saddle-point at $u_a = (s_{12})^{-1/2}$.
We emphasize that the $v$-integral cannot be approximated by the saddle-point method, we will compute it separately later.
In summary, we transform the $(a_1, a_2)$-integral as follows: ($f$ = integrand)
\begin{equation}
    \int_0^\infty \dd a_1 \int_0^\infty \dd a_2 \, f(a_1, a_2)
    = \int_0^\infty \dd u_a \int_0^\infty \dd v \, f\p2{\eta^{1/3} u_a v, \eta^{-1/3} \frac{u_a}{v}} 2 \frac{u_a}{v}.
\end{equation}
For the $b$-integral we introduce the new coordinates $u_b \in U(1)$ and $w \in U(1)$ and write $\wt{b}_1 := u_b w$, $\wt{b}_2 := \frac{u_b}{w}$.
Transforming the integrals via $\oint \dd b \, f(b) = c \oint \dd b \, f(c b)$ and
\begin{eqsplit*}
    \oint \dd \wt{b}_1 \oint \dd \wt{b}_2 \, f(\wt{b}_1, \wt{b}_2)
    & = \int_{-\pi}^\pi \dd t \int_{-\pi}^\pi \dd s \, f(\ee^{\ii t}, \ee^{\ii s}) \ii \ee^{\ii t} \, \ii \ee^{\ii s}
    = \frac{1}{2} \int_{-\pi}^\pi \dd x \int_{-\pi}^\pi \dd y \, f(\ee^{\ii (x + y)}, \ee^{\ii (x - y)}) \ii^2 \ee^{2 \ii x} 2
    \\
    & = \oint \dd u_b \oint \dd w \, f\p2{u_b w, \frac{u_b}{w}} \frac{u_b}{w}
\end{eqsplit*}
shows that
\begin{equation}
    \oint \dd b_1 \oint \dd b_2 \, f(b_1, b_2)
    = \oint \dd u_b \oint \dd w \, f\p2{\eta^{1/3} u_b w, \eta^{-1/3} \frac{u_b}{w}} \frac{u_b}{w}
\end{equation}
The $u_b$-integral has two saddle-points: $u_b = \pm (s_{12})^{-1/2}$.
In summary, we get
\begin{eqsplit}
    \E\s*{\Tr\p*{\frac{1}{H - \eta \zeta}}}
    & = -\frac{2 \ii N^3}{(2 \pi)^2} \int_0^\infty \dd v \oint \dd w
    \int_0^\infty \dd u_a \oint \dd u_b \, \ee^{
        -N I\p*{\eta^{1/3} u_a v, \eta^{-1/3} \frac{u_a}{v}}
        + N I\p*{\eta^{1/3} u_b w, \eta^{-1/3} \frac{u_b}{w}}
    } \\
    & \qquad \times \det\s*{ S + \diag\p*{\frac{1}{\eta^{2/3} u_a v u_b w}, \frac{\eta^{2/3} v w}{u_a u_b}} } \p2{\eta^{1/3} u_a v + \eta^{-1/3} \frac{u_a}{v}} \frac{u_a u_b}{v w}.
\end{eqsplit}

We carry out one more step to simplify this integral representation.
First, we perform the transformations $u_a \to \frac{1}{\sqrt{s_{12}}} u_a$, $u_b \to \frac{1}{\sqrt{s_{12}}} u_b$, $v \to \sqrt{\frac{s_{12}}{s_{11}}} \, v$ and $w \to \sqrt{\frac{s_{12}}{s_{11}}} \, w$.
Then we introduce the functions
\begin{align}
    J(u) & := u_a^2 - 2 \log(u_a) - u_b^2 + 2 \log(u_b),
    \\
    L(u; v, w) & := N \eta^{2/3} \p*{
        \frac{1}{2} (u_a^2 v^2 - u_b^2 w^2)
        - \ii \zeta \p2{\frac{u_a}{v} - \frac{u_b}{w}}
    }
    + N \eta^{4/3} \frac{s_{12}}{s_{11}} (-\ii \zeta) (u_a v - u_b w),
    \\
    F(u; v, w) & := \ee^{-L(u; v, w)} \p2{\eta^{1/3} \frac{s_{12}}{s_{11}} \, u_a v + \eta^{-1/3} \frac{u_a}{v}} \frac{u_a u_b}{v w},
    \\
    D(u; v, w) & := \frac{1}{s_{12}^2} \det\s*{ S + \diag\p*{\frac{s_{11}}{\eta^{2/3} u_a v u_b w}, \frac{s_{12}^2}{s_{11}} \frac{\eta^{2/3} v w}{u_a u_b}} },
\end{align}
where $u := (u_a, u_b) \in (0, \infty) \times U(1)$ for brevity.
Now we can write
\begin{eqsplit}
    \E\s*{\Tr\p*{\frac{1}{H - \eta \frac{s_{12}}{\sqrt{s_{11}}} \zeta}}}
    & = -\frac{2 \ii N^3}{(2 \pi)^2} \int_0^\infty \dd v \oint \dd w
    \int_0^\infty \dd u_a \oint \dd u_b \, \ee^{-N J(u)} \frac{s_{11}^{3/2}}{s_{12}^3} F(u; v, w) s_{12}^2 D(u; v, w) \frac{1}{s_{11}}.
\end{eqsplit}
As $\eta \to 0$ the function $J$ contains the leading order terms of $I(a) - I(b)$, while $L$ contains the lower order terms times $N$.
In particular, for $\eta \asymp N^{-3/2}$ we have $L(u; v, w) = O(1)$.
Therefore, the saddle-points $u^* = (u_a^*, u_b^*)$ of the $u$-integrals only depend on $J$.
We chose the transformation of $u_a$ and $u_b$ such that the saddle-points are given by $u_a^* = 1$ and $u_b^* = \pm 1$.
Then we chose the rescaling of $v$, $w$ and $\zeta$ such that the leading order of $L$ is independent of $s_{11}$ and $s_{12}$.
Finally, we defined $D$ as the determinant and $F$ as the remaining terms in the integrand, choosing appropriate prefactors to make their leading order independent of $s_{11}$ and $s_{12}$.
For the $D$-term it will be useful to expand the determinant in $\eta$:
\begin{equation}
    D(u; v, w) = D^{(0)}(u) + \eta^{2/3} D^{(1)}(u; v, w), \qquad \text{where }
    D^{(0)}(u) := \frac{1}{u_a^2 u_b^2} - 1, \quad
    D^{(1)}(u; v, w) := \frac{v w}{u_a u_b}.
\end{equation}
Note that the Jacobian of the $(u, v, w)$-transforms is $s_{11}^{-1}$.
The product of the $s_{11}$ and $s_{12}$ dependent factors leads to the rescaling on the left-hand side of \eqref{eq: K=2 expectation limit}.
In particular, the right-hand side of \eqref{eq: K=2 expectation limit} is universal, provided that the sub-leading terms of $F$ and $L$ (which still depend on $s_{11}$ and $s_{12}$) do not contribute to the leading order of the integral.
As we will see in the following calculations, the integral over the leading order of $F$ and $L$ is non-zero, confirming this assumption.

\subsubsection{Approximation of the \texorpdfstring{$u$}{u}-integral}
The critical points and saddle directions of $J(u_a, u_b)$ are
\[
    u^* = \begin{cases}
        (1, +1), \\
        (1, -1),
    \end{cases}
    \qquad \text{and} \qquad
    \mu = \begin{cases}
        (1, +\ii), \\
        (1, -\ii),
    \end{cases}
    = (1, \ii u_b^*)
\]
with second derivatives (note that we substitute $u_a^* = 1$ and $(u_b^*)^2 = 1$ where possible)
\[
    \mu_a^2 (\partial_{u_a}^2 J)(u^*) = 4, \quad
    \mu_b^2 (\partial_{u_b}^2 J)(u^*) = 4.
\]
We also have $\ee^{-N J(u^*)} = 1$ for $N \in \N$ and the third derivatives
\[
    (\partial_{u_a}^3 J)(u^*) = -4, \quad
    (\partial_{u_b}^3 J)(u^*) = \frac{4}{u_b^*}.
\]
For the functions $F$ and $L$ we will only use the leading terms as $\eta \to 0$
\begin{align*}
    F(u_a, u_b; v, w) & = \eta^{-1/3} \ee^{-L(u_a, u_b; v, w)} \frac{u_a^2 u_b}{v^2 w} + O(\eta^{1/3}),
    \\
    L(u_a, u_b; v, w) & = N \eta^{2/3} \p2{
        \frac{1}{2} (u_a^2 v^2 - u_b^2 w^2)
        - \ii \zeta \p2{\frac{u_a}{v} - \frac{u_b}{w}}
    } + O(N \eta^{4/3})
\end{align*}
and the derivatives
\begin{align*}
    (\partial_{u_a} F)(u^*; v, w) & = \eta^{-1/3} \ee^{-L(u^*; v, w)} \frac{1}{v^2 w} \s1{-(\partial_{u_a} L)(u^*; v, w) + 2} u_b^* + O(\eta^{1/3}),
    \\
    (\partial_{u_b} F)(u^*; v, w) & = \eta^{-1/3} \ee^{-L(u^*; v, w)} \frac{1}{v^2 w} \s1{(\partial_{u_b} L)(u^*; v, w) \, u_b^* + 1} u_b^* + O(\eta^{1/3}),
    \\
    (\partial_{u_a} L)(u^*; v, w) & = N \eta^{2/3} \p2{v^2 - \ii \zeta \frac{1}{v}} + O(N \eta^{4/3}),
    \\
    (\partial_{u_b} L)(u^*; v, w) & = N \eta^{2/3} \p2{w^2 u_b^* - \ii \zeta \frac{1}{w}} + O(N \eta^{4/3}).
\end{align*}

For $D^{(0)}$ and $D^{(1)}$ we compare with the formulas in Appendix~\ref{app: Saddle-Point Method} to determine which derivatives we need.
As in the appendix, we use a subscript to indicate the order of the total derivatives in $u$, evaluated at $u = u^*$.
Then for $D^{(0)}$ we have $D^{(0)}_0 = 0$ (since $(u_a^* u_b^*)^2 = 1$) and $D^{(0)}_1 \neq 0$.
Thus, the integral over $D^{(0)}$ is of order $\eta^{-1/3} O(N^{-1})$ (with the $\eta^{-1/3}$ coming from $F$).
For $D^{(1)}$ we have $D^{(1)}_0 \neq 0$, but since $N \eta^{2/3} \asymp 1$ by our choice of $\eta$, we see that the integral over $\eta^{2/3} D^{(1)}$ is also of order $\eta^{-1/3} O(N^{-1})$.
We conclude that the integrals over $D^{(0)}$ and $\eta^{2/3} D^{(1)}$ both contribute to the leading order of the saddle-point approximation.
The required derivatives for the $D^{(0)}$-integral are
\begin{gather*}
    (\partial_{u_a} D^{(0)})(u^*) = -2, \quad
    (\partial_{u_b} D^{(0)})(u^*) = \frac{2}{u_b^*},
    \\
    (\partial_{u_a}^2 D^{(0)})(u^*) = 6, \quad
    (\partial_{u_a} \partial_{u_b} D^{(0)})(u^*) = \frac{4}{u_b^*}, \quad
    (\partial_{u_b}^2 D^{(0)})(u^*) = 6.
\end{gather*}
For the $D^{(1)}$-integral we only need the value $D^{(1)}(u^*; v, w) = \frac{v w}{u_b^*}$.

Let us introduce the notation $f(N) \sim \phi(N)$ when $\lim_{N \to \infty} \frac{f(N)}{\phi(N)} = 1$ (asymptotic equality) and $f(N) = O(\phi(N))$ when $\limsup_{N\ \to \infty} \abs*{\frac{f(N)}{\phi(N)}} < \infty$.
Using the saddle-point method we obtain the following approximations
\begin{eqsplit*}
    \MoveEqLeft \int_0^\infty \dd u_a \oint \dd u_b \, \ee^{-N J} F D^{(0)}
    \\
    & = \sum_{u^*} \frac{\mu_a \mu_b}{\sqrt{N}^2} \frac{1}{N} \ee^{-N J_0} \s*{
        \int_{\R^2} \dd y \, \ee^{-J_2} \p*{F_1 D^{(0)}_1 + F_0 D^{(0)}_2 - F_0 D^{(0)}_1 J_3}
        + O\p*{\frac{1}{N}}
    }
    \\
    & \sim \frac{\ii \pi}{4 N^2} \eta^{-1/3} \frac{1}{v^2 w} \sum_{u_b^* = \pm 1} u_b^* \ee^{-L(1, u_b^*; v, w)} \s*{
        (\partial_{u_a} L)(1, u_b^*; v, w) \, u_b^* - 1 + (\partial_{u_b} L)(1, u_b^*; v, w) - u_b^*
    }
    \\
    & \sim \frac{\ii \pi}{4 N^2} \eta^{-1/3} \ee^{
        -N \eta^{2/3} \p*{\frac{1}{2} v^2 - \frac{\ii \zeta}{v}}
    } \sum_{u_b^* = \pm 1} \ee^{
        N \eta^{2/3} \p*{\frac{1}{2} w^2 - \frac{\ii \zeta}{w} u_b^*}
    } N \eta^{2/3}
    \\
    & \qquad \times \s*{
        \frac{1}{w} - \frac{\ii \zeta}{v^3 w} - \frac{1}{N \eta^{2/3} v^2 w} + \frac{w}{v^2} - \frac{\ii \zeta}{v^2 w^2} u_b^*
    }
\end{eqsplit*}
and
\begin{eqsplit*}
    \MoveEqLeft \int_0^\infty \dd u_a \oint \dd u_b \, \ee^{-N J} F D^{(1)}
    = \sum_{u^*} \frac{\mu_a \mu_b}{\sqrt{N}^2} \ee^{-N J_0} \s*{
        \int_{\R^2} \dd y \, \ee^{-J_2} F_0 D_0^{(1)}
        + O\p*{\frac{1}{N}}
    }
    \\
    & \sim \frac{\ii \pi}{2 N} \eta^{-1/3} \frac{1}{v^2 w} \sum_{u_b^* = \pm 1} u_b^* \ee^{-L(1, u_b^*; v, w)} v w
    \\
    & \sim \frac{\ii \pi}{2 N} \eta^{-1/3} \ee^{
        -N \eta^{2/3} \p*{\frac{1}{2} v^2 - \frac{\ii \zeta}{v}}
    } \sum_{u_b^* = \pm 1} \ee^{
        N \eta^{2/3} \p*{\frac{1}{2} w^2 - \frac{\ii \zeta}{w} u_b^*}
    } \frac{1}{v} u_b^*.
\end{eqsplit*}
For the expectation of the resolvent we arrive at
\begin{eqsplit*}
    \MoveEqLeft \frac{s_{12}}{\sqrt{s_{11}}} \E\s*{\Tr\p*{\frac{1}{H - \eta \frac{s_{12}}{\sqrt{s_{11}}} \zeta}}}
    = -\frac{\ii N^3}{2 \pi^2} \int_0^\infty \dd v \oint \dd w \int_0^\infty \dd u_a \oint \dd u_b \, \ee^{-N J} F \p*{D^{(0)} + \eta^{2/3} D^{(1)}}
    \\
    & \sim \frac{N^2}{8 \pi} \eta^{1/3} \int_{0}^{\infty} \dd v \, \ee^{
        -N \eta^{2/3} \p*{\frac{1}{2} v^2 - \frac{\ii \zeta}{v}}
    } \sum_{u_b^* = \pm 1} \oint \dd w \, \ee^{
        N \eta^{2/3} \p*{\frac{1}{2} w^2 - \frac{\ii \zeta}{w} u_b^*}
    }
    \\
    & \qquad \times \s*{
        \frac{1}{w} - \frac{\ii \zeta}{v^3 w} - \frac{1}{N \eta^{2/3} v^2 w} + \frac{w}{v^2} + \p*{- \frac{\ii \zeta}{v^2 w^2} + 2 \frac{1}{v}} u_b^*
    }.
\end{eqsplit*}

\subsubsection{Computation of the \texorpdfstring{$v$}{v}- and \texorpdfstring{$w$}{w}-integrals}
For the next steps, let us temporarily introduce the variables $t_1 := N \eta^{2/3} \frac{1}{2}$ and $t_2 := - N \eta^{2/3} \ii \zeta$.
To solve the $v$-integral we perform the substitution $v = \sqrt{x}$, then write the integrand in terms of Meijer $G$-functions and apply the $G$-function convolution formula:
\begin{align*}
    \int_0^\infty \dd v \, \ee^{-(t_1 v^2 + \frac{t_2}{v})} v^n
    & = \frac{1}{2} \int_0^\infty \dd x \, x^{\frac{n}{2} - \frac{1}{2}} \exp\p*{-t_1 x} \times \exp\p*{- \frac{t_2}{\sqrt{x}}}
    \\
    & = \frac{1}{2} \int_0^\infty \dd x \, \frac{1}{t_1^{(n-1)/2}} \MeijerG{1, 0}{0, 1}{-}{\frac{n - 1}{2}}{t_1 x} \times \frac{1}{\sqrt{\pi}} \MeijerG{0, 2}{2, 0}{\frac{1}{2}, 1}{-}{\frac{4 x}{t_2^2}}
    \\
    & = \frac{1}{2 \sqrt{\pi}} t_1^{- \frac{n - 1}{2}} \frac{t_2^2}{4} \MeijerG{3, 0}{0, 3}{-}{\frac{n - 1}{2}, -\frac{1}{2}, -1}{\frac{t_1}{4 / t_2^2}}.
\end{align*}
Using \cite[\DLMFeqref{16.19}{2}]{NIST:DLMF} we can further simplify the result and arrive at
\begin{equation} \label{eq: v-integral Meijer G}
    \int_0^\infty \dd v \, \ee^{-(t_1 v^2 + \frac{t_2}{v})} v^n
    = \frac{1}{2 \sqrt{\pi}} t_1^{-\frac{n + 1}{2}} \MeijerG{3, 0}{0, 3}{-}{0, \frac{1}{2}, \frac{n + 1}{2}}{\frac{t_1 t_2^2}{4}}, \qquad
    \text{for } \Re(t_1) > 0,\ \Re(t_2) > 0,\ n \in \Z.
\end{equation}
For the $w$-integral we insert the series representations of $\ee^{t_1 w^2}$ and $\ee^{t_2 u_b^* / w}$.
Then we use that $\oint \dd w \, w^n = 2 \pi \ii \delta_{n, -1}$ for $n \in \Z$ and obtain
\begin{eqsplit}
    \oint \dd w \, \ee^{t_1 w^2 + \frac{t_2}{w} u_b^*} w^r
    & = \sum_{k = 0}^\infty \frac{t_1^k}{k!} \sum_{l = 0}^\infty \frac{(t_2 u_b^*)^l}{l!} \oint \dd w \, w^{2 k - l + r}
    \\
    & = 2 \pi \ii (t_2 u_b^*)^{1 + r} \sum_{k = 0}^{\infty} \frac{(t_1 (t_2 u_b^*)^2)^k}{k! \, (2 k + 1 + r)!}, \qquad
    \text{for } t_1, t_2 \in \C,\ r \in \Z.
\end{eqsplit}
The infinite series can be written as a generalized hypergeometric function:
\begin{equation} \label{eq: w-sum Hypergeometric F}
    \sum_{k = 0}^{\infty} \frac{1}{(2 k + 1 + r)!} \frac{x^k}{k!}
    = \begin{dcases}
        \frac{x^{-(1+r)/2}}{(-(1+r)/2)!} \, \HypergeometricPFQ{0}{2}{-}{\frac{1}{2}, \frac{1 - r}{2}}{\frac{x}{4}}, & r < -1,\ r \text{ odd},
        \\
        \frac{x^{-r/2}}{(-r/2)!} \, \HypergeometricPFQ{0}{2}{-}{\frac{3}{2}, \frac{2 - r}{2}}{\frac{x}{4}}, & r < -1,\ r \text{ even},
        \\
        \frac{1}{(1 + r)!} \, \HypergeometricPFQ{0}{2}{-}{\frac{2 + r}{2}, \frac{3 + r}{2}}{\frac{x}{4}}, & r \geq -1.
    \end{dcases}
\end{equation}
Note that $x = t_1 (t_2 u_b^*)^2 = t_1 t_2^2$ is independent of $u_b^* = \pm 1$.
Therefore
\[
    \sum_{u_b^* = \pm 1} \oint \dd w \, \ee^{t_1 w^2 + \frac{t_2}{w} u_b^*} w^r (u_b^*)^l
    = \p*{\oint \dd w \, \ee^{t_1 w^2 + \frac{t_2}{w}} w^r} \times \sum_{u_b^* = \pm 1} (u_b^*)^{1 + r + l}
\]
and the sum is either 0 (if $r + l$ is even) or 2 for $u_b^* = 1$ (if $r + l$ is odd).
After solving the integrals, the expectation of the resolvent is
\begin{eqsplit*}
    \MoveEqLeft \frac{s_{12}}{\sqrt{s_{11}}} \E\s*{\Tr\p*{\frac{1}{H - \eta \frac{s_{12}}{\sqrt{s_{11}}} \zeta}}}
    \sim \frac{N^2}{8 \pi} \eta^{1/3} \frac{2 \pi \ii}{2 \sqrt{\pi}} 2
    \\
    & \times \bigg[
        \sqrt{2} (N \eta^{2/3})^{-1/2} \,
        \MeijerG{3, 0}{0, 3}{-}{0, \frac{1}{2}, \frac{1}{2}}{\frac{N^3 \eta^2}{8} (-\ii \zeta)^2}
        \HypergeometricPFQ{0}{2}{-}{\frac{1}{2}, 1}{\frac{N^3 \eta^2}{8} (-\ii \zeta)^2}
    \\
    & \hphantom{\times \bigg[}
        + \frac{1}{2} (N \eta^{2/3}) (-\ii \zeta) \,
        \MeijerG{3, 0}{0, 3}{-}{0, \frac{1}{2}, -1}{\frac{N^3 \eta^2}{8} (-\ii \zeta)^2}
        \HypergeometricPFQ{0}{2}{-}{\frac{1}{2}, 1}{\frac{N^3 \eta^2}{8} (-\ii \zeta)^2}
    \\
    & \hphantom{\times \bigg[}
        - \sqrt{\frac{1}{2}} (N \eta^{2/3})^{-1/2} \,
        \MeijerG{3, 0}{0, 3}{-}{0, \frac{1}{2}, -\frac{1}{2}}{\frac{N^3 \eta^2}{8} (-\ii \zeta)^2}
        \HypergeometricPFQ{0}{2}{-}{\frac{1}{2}, 1}{\frac{N^3 \eta^2}{8} (-\ii \zeta)^2}
    \\
    & \hphantom{\times \bigg[}
        + 2 \p*{\frac{1}{2}}^{3/2} (N \eta^{2/3})^{5/2} (-\ii \zeta)^2 \,
        \MeijerG{3, 0}{0, 3}{-}{0, \frac{1}{2}, -\frac{1}{2}}{\frac{N^3 \eta^2}{8} (-\ii \zeta)^2}
        \HypergeometricPFQ{0}{2}{-}{\frac{3}{2}, 2}{\frac{N^3 \eta^2}{8} (-\ii \zeta)^2}
    \\
    & \hphantom{\times \bigg[}
        + 2 (N \eta^{2/3}) (-\ii \zeta) \,
        \MeijerG{3, 0}{0, 3}{-}{0, \frac{1}{2}, 0}{\frac{N^3 \eta^2}{8} (-\ii \zeta)^2}
        \HypergeometricPFQ{0}{2}{-}{1, \frac{3}{2}}{\frac{N^3 \eta^2}{8} (-\ii \zeta)^2}
    \bigg].
\end{eqsplit*}
To further simplify this expression, we use \cite[\DLMFeqref{16.19}{2}]{NIST:DLMF} to absorb any $-\ii \zeta$ factors into the Meijer $G$-functions.
Then we use the identity
\begin{equation*}
    \MeijerG{3, 0}{0, 3}{-}{\frac{1}{2}, 1, -\frac{1}{2}}{z} - \frac{1}{2} \MeijerG{3, 0}{0, 3}{-}{0, \frac{1}{2}, -\frac{1}{2}}{z}
    = \MeijerG{3, 0}{0, 3}{-}{0, \frac{1}{2}, \frac{1}{2}}{z}
\end{equation*}
for the second and third term and combine the result with the first term.
This identity follows from the integral representation \cite[\DLMFeqref{16.17}{1}]{NIST:DLMF}
\[
    \MeijerG{3, 0}{0, 3}{-}{b_1, b_2, b_3}{z} = \frac{1}{2 \pi \ii} \int_L \dd s \, \Gamma(b_1 - s) \Gamma(b_2 - s) \Gamma(b_3 - s) z^s
\]
and by merging the integrands as follows:
\[
    \Gamma\p*{-\frac{1}{2} - s} \Gamma\p*{\frac{1}{2} - s} \Gamma(1 - s) - \frac{1}{2} \Gamma\p*{-\frac{1}{2} - s} \Gamma\p*{\frac{1}{2} - s} \Gamma(0 - s)
    = \Gamma(0 - s) \Gamma\p*{\frac{1}{2} - s} \Gamma\p*{\frac{1}{2} - s}.
\]

The final result is
\begin{eqsplit*}
    \frac{s_{12}}{\sqrt{s_{11}}} \E\s*{\Tr\p*{\frac{1}{H - \eta \frac{s_{12}}{\sqrt{s_{11}}} \zeta}}}
    & \sim \frac{\ii \sqrt{2}}{\sqrt{\pi}} N^{3/2} \bigg[
        \frac{1}{2}
        \MeijerG{3, 0}{0, 3}{-}{0, \frac{1}{2}, \frac{1}{2}}{\frac{N^3 \eta^2}{8} (-\ii \zeta)^2}
        \HypergeometricPFQ{0}{2}{-}{\frac{1}{2}, 1}{\frac{N^3 \eta^2}{8} (-\ii \zeta)^2}
    \\
    & \hphantom{\sim \frac{\ii \sqrt{2}}{\sqrt{\pi}} N^{3/2} \bigg[}
        +
        \MeijerG{3, 0}{0, 3}{-}{\frac{1}{2}, 1, \frac{3}{2}}{\frac{N^3 \eta^2}{8} (-\ii \zeta)^2}
        \HypergeometricPFQ{0}{2}{-}{\frac{3}{2}, 2}{\frac{N^3 \eta^2}{8} (-\ii \zeta)^2}
    \\
    & \hphantom{\sim \frac{\ii \sqrt{2}}{\sqrt{\pi}} N^{3/2} \bigg[}
        +
        \MeijerG{3, 0}{0, 3}{-}{\frac{1}{2}, \frac{1}{2}, 1}{\frac{N^3 \eta^2}{8} (-\ii \zeta)^2}
        \HypergeometricPFQ{0}{2}{-}{1, \frac{3}{2}}{\frac{N^3 \eta^2}{8} (-\ii \zeta)^2}
    \bigg].
\end{eqsplit*}
Observe that the first two terms originate from $D^{(0)}$ and the last term from $D^{(1)}$.

\subsection{\texorpdfstring{$3 \times 3$}{3x3} case}
Let $S$ be of the form
\[
    S = \begin{pmatrix}
        s_{11} & s_{12} & s_{13} \\
        s_{12} & s_{22} & 0 \\
        s_{13} & 0 & 0
    \end{pmatrix}, \quad s_{11} \geq 0,\ s_{12}, s_{22}, s_{13} > 0
\]
Then the saddle-point equations $\nabla I(a) = 0$ are equivalent to
\begin{equation} \label{eq: K=3 saddle-point equation}
    \frac{1}{a_1} = s_{11} a_1 + s_{12} a_2 + s_{13} a_3 - \ii z, \qquad
    \frac{1}{a_2} = s_{12} a_1 + s_{22} a_2 - \ii z, \qquad
    \frac{1}{a_3} = s_{13} a_1 - \ii z.
\end{equation}
For $z = \eta \zeta$ with $\eta \to 0$ we have that $\abs{a_1} \asymp \eta^{1/2}$, $\abs{a_2} \asymp 1$ and $\abs{a_3} \asymp \eta^{-1/2}$.
We refer to Section~\ref{sec: Weak Non-Chirality Limit for K = 3} for a procedure to determine these rates.

\subsubsection{Transformation of the integrals}
We introduce $\wt{a}_1 := \eta^{-1/2} a_1$ and $\wt{a}_3 := \eta^{1/2} a_3$, then
\[
    I(a) = s_{13} \wt{a}_1 \wt{a}_3 + \frac{s_{22}}{2} a_2^2 - \log(\wt{a}_1) - \log(\wt{a}_2) - \log(\wt{a}_3) + O\p*{\eta^{1/2}} \qquad \text{for } \eta \to 0.
\]
Now the $a_2$-integral has a saddle-point at $a_2 = s_{22}^{-1/2}$ and the $(\wt{a}_1, \wt{a}_3)$-integral has a saddle-line at $\wt{a}_3 = (s_{13} \wt{a}_1)^{-1}$.
Analogous to the $2 \times 2$-case, we transform the saddle-line into a saddle-point by introducing new coordinates $u_a$ and $v$.
To maintain uniform notation, we also replace $a_2$ with a $u$-variable.
The complete transformation is
\[
    (a_1, a_2, a_3) = \p*{\eta^{1/2} u_{a1} v, u_{a2}, \eta^{-1/2} \frac{u_{a1}}{v}}, \qquad
    \text{where } u_a = (u_{a1}, u_{a2}) \in (0, \infty)^2 \text{ and } v > 0.
\]
Similarly for the $b$-integral we write
\[
    (b_1, b_2, b_3) = \p*{\eta^{1/2} u_{b1} w, u_{b2}, \eta^{-1/2} \frac{u_{b1}}{w}}, \qquad
    \text{where } u_b = (u_{b1}, u_{b2}) \in U(1)^2 \text{ and } w \in U(1).
\]
For the resulting transformation of the $a$- and $b$-integrals, see the $2 \times 2$-case.
In summary, we get
\begin{eqsplit}
    \E\s*{\Tr\p*{\frac{1}{H - \eta \zeta}}}
    & = -\frac{2 N^4}{(2 \pi)^3} \int_0^\infty \dd v \oint \dd w
        \int_{(0, \infty)^2} \dd u_a \oint^2 \dd u_b \\
    & \qquad \times \ee^{
            -N I\p*{\eta^{1/2} u_{a1} v, u_{a2}, \eta^{-1/2} \frac{u_{a1}}{v}}
            + N I\p*{\eta^{1/2} u_{b1} w, u_{b2}, \eta^{-1/2} \frac{u_{b1}}{w}}
        } \\
    & \qquad \times \det\s*{ S + \diag\p*{\frac{1}{\eta u_{a1} v u_{b1} w}, \frac{1}{u_{a2} u_{b2}}, \frac{\eta v w}{u_{a1} u_{b1}}} } \\
    & \qquad \times \p2{\eta^{1/2} u_{a1} v + u_{a2} + \eta^{-1/2} \frac{u_{a1}}{v}} \frac{u_{a1} u_{b1}}{v w}.
\end{eqsplit}

To simplify the integral representation, we perform the transformations
\begin{gather*}
    u_{a1} \to \frac{1}{\sqrt{s_{13}}} u_{a1}, \quad
    u_{a2} \to \frac{1}{\sqrt{s_{22}}} u_{a2}, \quad
    u_{b1} \to \frac{1}{\sqrt{s_{13}}} u_{b1}, \quad
    u_{b2} \to \frac{1}{\sqrt{s_{22}}} u_{b2}, \\
    v \to \frac{\sqrt{s_{13} s_{22}}}{s_{12}} v, \quad
    w \to \frac{\sqrt{s_{13} s_{22}}}{s_{12}} w
\end{gather*}
and introduce the functions
\begin{align*}
    J(u) & := u_{a1}^2 - 2 \log(u_{a1}) + \frac{1}{2} u_{a2}^2 - \log(u_{a2}) - u_{b1}^2 + 2 \log(u_{b1}) - \frac{1}{2} u_{b2}^2 + \log(u_{b2}),
    \\
    L(u; v, w) & := N \eta^{1/2} \p2{
        (u_{a1} u_{a2} v - u_{b1} u_{b2} w)
        - \ii \zeta \p2{\frac{u_{a1}}{v} - \frac{u_{b1}}{w}}
    }
    \\
    & \quad + N \eta \p2{
        \frac{s_{11} s_{22}}{2 s_{12}^2} (u_{a1}^2 v^2 - u_{b1}^2 w^2)
        + \frac{s_{13}}{s_{12}} (- \ii \zeta) (u_{a2} - u_{b2})
    }
    + N \eta^{3/2} \frac{\sqrt{s_{22}}}{s_{12}} (- \ii \zeta) (u_{a1} v - u_{b1} w),
    \\
    F(u; v, w) & := \ee^{-L(u; v, w)} \p2{\eta^{1/2} \frac{s_{13} s_{22}}{s_{12}^2} u_{a1} v + \frac{s_{13}}{s_{12}} u_{a2} + \eta^{-1/2} \frac{u_{a1}}{v}} \frac{u_{a1} u_{b1}}{v w},
    \\
    D(u; v, w) & := \frac{1}{s_{13}^2 s_{22}} \det\s*{ S + \diag\p*{\frac{s_{12}^2}{s_{22}} \frac{1}{\eta u_{a1} v u_{b1} w}, \frac{s_{22}}{u_{a2} u_{b2}}, \frac{s_{13}^2 s_{22}}{s_{12}^2} \frac{\eta v w}{u_{a1} u_{b1}}} }.
\end{align*}
where $u := (u_a, u_b) \in (0, \infty)^2 \times U(1)^2$ for brevity.
Now we can write
\begin{eqsplit}
    \E\s*{\Tr\p*{\frac{1}{H - \eta \frac{s_{13} \sqrt{s_{22}}}{s_{12}} \zeta}}}
    & = -\frac{N^4}{4 \pi^3}
        \int_0^\infty \dd v \oint \dd w
        \int_{(0, \infty)^2} \dd u_a \oint^2 \dd u_b \, \ee^{-N J(u)}
    \\
    & \qquad \times \frac{s_{12}^3}{s_{13}^3 s_{22}^{3/2}} F(u; v, w) s_{13}^2 s_{22} D(u; v, w) \frac{1}{s_{12}^2}.
\end{eqsplit}
Analogous to the $K = 2$ case, the function $J$ contains the leading order of $I(a) - I(b)$, while $L$ contains the lower order terms times $N$.
In particular, for $\eta \asymp N^{-2}$ we have $L(u; v, w) = O(1)$.
We chose the transformations of $u_a$ and $u_b$ such that $J(u)$ is independent of $s_{13}$ and $s_{22}$, then we transformed $v$, $w$ and $\zeta$ such that the leading order of $L(u; v, w)$ is independent of $s_{12}$, $s_{13}$ and $s_{22}$.
Finally, we defined $F$ and $D$ such that their leading orders are also independent of the parameters in $S$.
The $D$-term can be expanded in $\eta$ as follows: $D(u; v, w) = D^{(0)}(u) + \eta D^{(1)}(u; v, w)$ with
\begin{equation}
    D^{(0)}(u) := \p*{\frac{1}{u_{a1}^2 u_{b1}^2} - 1} \p*{\frac{1}{u_{a2} u_{b2}} + 1}, \quad
    D^{(1)}(u; v, w) := \frac{v w}{u_{a1} u_{b1}} \p*{\frac{s_{11} s_{22}}{s_{12}^2} \p*{1 + \frac{1}{u_{a2} u_{b2}}} - 1}.
\end{equation}
We emphasize the subleading orders of $L$, $F$ and $D$ still depend on the parameters in $S$.
The calculations below show that the integral over the leading orders of $F$ and $L$ is not zero and that the subleading order of $D$ does not contribute to the leading order of the saddle-point approximation (unlike the $K = 2$ case).
Thus, all transformations taken together lead to the rescaling on the left-hand side of \eqref{eq: K=3 expectation limit}.
In particular, the right-hand side of \ref{eq: K=3 expectation limit} is universal.

\subsubsection{Approximation of the \texorpdfstring{$u$}{u}-integral}
The critical points and saddle directions of $J(u)$ are
\[
    u^* = \begin{cases}
        (1, 1, +1, +1), \\
        (1, 1, +1, -1), \\
        (1, 1, -1, +1), \\
        (1, 1, -1, -1),
    \end{cases}
    \qquad \text{and} \qquad
    \mu = \begin{cases}
        (1, 1, +\ii, +\ii), \\
        (1, 1, +\ii, -\ii), \\
        (1, 1, -\ii, +\ii), \\
        (1, 1, -\ii, -\ii), \\
    \end{cases}
    = (1, 1, \ii u_{b1}^*, \ii u_{b2}^*)
\]
with second derivatives (we substitute $u_{a1}^* = 1$, $u_{a2}^* = 1$, $(u_{b1}^*)^2 = 1$ and $(u_{b2}^*)^2 = 1$ where possible)
\[
    \mu_{a1}^2 (\partial_{u_{a1}}^2 J)(u^*) = 4, \quad
    \mu_{a2}^2 (\partial_{u_{a2}}^2 J)(u^*) = 2, \quad
    \mu_{b1}^2 (\partial_{u_{b1}}^2 J)(u^*) = 4, \quad
    \mu_{b2}^2 (\partial_{u_{b2}}^2 J)(u^*) = 2.
\]
We also need $\ee^{-N J(u^*)} = 1$ for $N \in \N$ and the third derivatives
\[
    (\partial_{u_{a1}}^3 J)(u^*) = -4, \quad
    (\partial_{u_{a2}}^3 J)(u^*) = -2, \quad
    (\partial_{u_{b1}}^3 J)(u^*) = \frac{4}{u_{b1}^*}, \quad
    (\partial_{u_{b2}}^3 J)(u^*) = \frac{2}{u_{b2}^*}.
\]
For the functions $F$ and $L$ we only require the $\eta \to 0$ leading terms
\begin{align*}
    F(u; v, w) & = \eta^{-1/2} \ee^{-L(u; v, w)} \frac{u_{a1}^2 u_{b1}}{v^2 w} + O(1),
    \\
    L(u; v, w) & = N \eta^{1/2} \p2{
        u_{a1} u_{a2} v - u_{b1} u_{b2} w
        - \ii \zeta \p2{\frac{u_{a1}}{v} - \frac{u_{b1}}{w}}
    } + O(N \eta)
\end{align*}
and the derivatives
\begin{align*}
    (\partial_{u_{a1}} F)(u^*; v, w) & = \eta^{-1/2} \frac{1}{v^2 w} \ee^{-L(u^*; v, w)} \s*{-(\partial_{u_{a1}} L)(u^*; v, w) + 2} u_{b1}^* + O(1),
    \\
    (\partial_{u_{a2}} F)(u^*; v, w) & = \eta^{-1/2} \frac{1}{v^2 w} \ee^{-L(u^*; v, w)} \s*{-(\partial_{u_{a2}} L)(u^*; v, w)} u_{b1}^* + O(1),
    \\
    (\partial_{u_{b1}} F)(u^*; v, w) & = \eta^{-1/2} \frac{1}{v^2 w} \ee^{-L(u^*; v, w)} \s*{1 + (\partial_{u_{b1}} L)(u^*; v, w) u_{b1}^*} + O(1),
    \\
    (\partial_{u_{b2}} F)(u^*; v, w) & = \eta^{-1/2} \frac{1}{v^2 w} \ee^{-L(u^*; v, w)} (\partial_{u_{b2}} L)(u^*; v, w) u_{b1}^* + O(1),
\end{align*}
\[
    (\partial_{u_{a1}} L)(u^*; v, w) \sim N \eta^{1/2} \p*{v u_{a2}^* - \frac{\ii \zeta}{v}} + O(N \eta), \qquad
    (\partial_{u_{a2}} L)(u^*; v, w) \sim N \eta^{1/2} v u_{a1}^* + O(N \eta).
\]

For $D^{(0)}$ and $D^{(1)}$ we compare with the formulas in Appendix~\ref{app: Saddle-Point Method} to determine which derivatives we need.
In the following we use a subscript to indicate the order of the derivative in $u$, evaluated at $u = u^*$.
Then for $D^{(0)}$ we have $D^{(0)}_0 = 0$ (since $(u_{a1}^* u_{b1}^*)^2 = 1$) and $D^{(0)}_1 \neq 0$ (we need at least one derivative in $u_{a1}$ or $u_{b1}$ to get a non-zero value).
Thus, the integral over $D^{(0)}$ is of order $\eta^{-1/2} O(N^{-1})$ (with the $\eta^{-1/2}$ factor coming from $F$).
Since $N \eta^{1/2} \asymp 1$, the order of the integral over $\eta D^{(1)}$ cannot be higher than $\eta^{-1/2} O(\eta) = \eta^{-1/2} O(N^{-2})$.
We conclude that the leading order of the saddle-point approximation is given by the integral over $D^{(0)}$ and that $D^{(1)}$ only contributes to sub-leading orders that we do not compute here.
The required derivatives for the $D^{(0)}$-integral are
\begin{eqsplit*}
    (\partial_{u_{a1}} D^{(0)})(u^*) & = -2 \p*{1 + \frac{1}{u_{b2}^*}}, \quad
    (\partial_{u_{a2}} D^{(0)})(u^*) = 0, \\
    (\partial_{u_{b1}} D^{(0)})(u^*) & = -2 \p*{1 + \frac{1}{u_{b2}^*}} u_{b1}^*, \quad
    (\partial_{u_{b2}} D^{(0)})(u^*) = 0.
\end{eqsplit*}

We obtain the approximation
\begin{eqsplit*}
    \MoveEqLeft \int_{(0^\infty)^2} \dd u_a \oint^2 \dd u_b \, \ee^{-N J} F D^{(0)} \\
    & = \sum_{u^*} \frac{\mu_{a1} \mu_{a2} \mu_{b1} \mu_{b2}}{\sqrt{N}^4} \frac{1}{N} \ee^{-N J_0} \s*{
        \int_{\R^4} \dd y \, \ee^{-J_2} \p*{F_1 D^{(0)}_1 + F_0 D^{(0)}_2 - F_0 D^{(0)}_1 J_3}
        + O\p*{\frac{1}{N}}
    }
    \\
    & \sim -\frac{\pi^2}{4 N^3} \sum_{u_{b1}^* = \pm 1, u_{b2}^* = \pm 1} u_{b1}^* u_{b2}^* \p*{1 + \frac{1}{u_{b2}^*}} \s*{
        -(\partial_{u_{a1}} F)(u^*) + (\partial_{u_{b1}} F)(u^*) u_{b1}^*
    } \Big\rvert_{u_{a1}^* = u_{a2}^* = 1}
    \\
    & \sim -\frac{2 \pi^2}{4 N^3} \eta^{-1/2} \frac{1}{v^2 w} \sum_{u_{b1}^* = \pm 1} \ee^{-L(u^*; v, w)} \big[
        (\partial_{u_{a1}} L)(u^*; v, w) (u_{b1}^*)^2 - 2 (u_{b1}^*)^2
        + (\partial_{u_{b1}} L)(u^*; v, w) u_{b1}^* + 1
    \big] \Big\rvert_{\begin{subarray}{l}u_{a1}^* = u_{a2}^* = 1, \\ u_{b2}^* = 1\end{subarray}}
    \\
    & \sim -\frac{\pi^2}{2 N^3} \eta^{-1/2} \ee^{-N \eta^{1/2}(v - \ii \zeta \frac{1}{v})} \frac{1}{v^2 w} \sum_{u_{b1}^* = \pm 1} \ee^{N \eta^{1/2}(w - \ii \zeta \frac{1}{w}) u_{b1}^*} N \eta^{1/2}
    \\
    & \qquad
    \times \s*{
        v + (-\ii \zeta) \frac{1}{v} - \frac{1}{N \eta^{1/2}} + \p*{w - \ii \zeta \frac{1}{w}} u_{b1}^*
    }.
\end{eqsplit*}
Observe that for $u_{b2}^* = -1$ we have $1 + \frac{1}{u_{b2}^*} = 0$ and for $u_{b1}^* = \pm 1$ we have $(u_{b1}^*)^2 = 1$.

For the expectation of the resolvent we arrive at
\begin{eqsplit*}
    \frac{s_{12}}{s_{13} \sqrt{s_{22}}} \E\s*{\Tr\p*{\frac{1}{H - \eta \frac{s_{13} \sqrt{s_{22}}}{s_{12}} \zeta}}}
    & \sim \frac{N^2}{4 \pi} \int_{0}^{\infty} \dd v \, \ee^{-N \eta^{1/2} \p*{v - \frac{\ii \zeta}{v}}} \oint \dd w \sum_{u_{b1}^* = \pm 1} \ee^{N \eta^{1/2} \p*{w - \frac{\ii \zeta}{w}} u_{b1}^*}
    \\
    & \quad \times \s*{
        \frac{1}{v w} - \frac{\ii \zeta}{v^3 w} - \frac{1}{N \eta^{1/2} v^2 w}
        + \p*{1 - \frac{\ii \zeta}{w^2}} u_{b1}^*
    }
\end{eqsplit*}

\subsubsection{Computation of the \texorpdfstring{$v$}{v}- and \texorpdfstring{$w$}{w}-integrals}
For the next steps, we introduce the variables $t_1 := N \eta^{1/2}$ and $t_2 := N \eta^{1/2} (-\ii \zeta)$.
The $v$-integral is given in \cite[(3.471.9)]{gradshteyn2000Table}.
To solve it by hand we can follow the same steps as in the $K = 2$ case: we write the integrand in terms of Meijer $G$-functions and then apply the $G$-function convolution theorem
\begin{align*}
    \int_0^\infty \dd v \, \ee^{-(t_1 v + \frac{t_2}{v})} v^n
    & = \int_0^\infty \dd v \, \frac{1}{t_1^n} \MeijerG{1, 0}{0, 1}{-}{n}{t_1 v} \times \MeijerG{0, 1}{1, 0}{1}{-}{\frac{v}{t_2}}
    = t_1^{- n} t_2 \MeijerG{2, 0}{0, 2}{-}{n, -1}{\frac{t_1}{1 / t_2}}.
\end{align*}
Using the identity $\MeijerG{2, 0}{0, 2}{-}{b, c}{z} = 2 z^{(b + c) / 2} K_{b - c}(2 \sqrt{z})$ from \cite[(9.34.3)]{gradshteyn2000Table}, where $K_\nu$ is the Bessel-$K$ function, we arrive at
\begin{equation} \label{eq: v-integral Bessel K}
    \int_0^\infty \dd v \, \ee^{-(t_1 v + \frac{t_2}{v})} v^n
    = 2 \p*{\frac{t_2}{t_1}}^{\frac{n + 1}{2}} K_{n + 1}\p*{2 \sqrt{t_1 t_2}}, \qquad
    \text{for } \Re(t_1) > 0,\ \Re(t_2) > 0,\ n \in \Z.
\end{equation}
For the $w$-integral we insert the series representations of $\ee^{t_1 w^2 u_{b1}^*}$ and $\ee^{t_2 u_{b1}^* / w}$.
Then we use that $\oint \dd w \, w^n = 2 \pi \ii \delta_{n, -1}$ for $n \in \Z$ and obtain
\begin{align*}
    \oint \dd w \, \ee^{(t_1 w + \frac{t_2}{w}) u_{b1}^*} w^r
    = \sum_{k = 0}^\infty \frac{t_1^k}{k!} \sum_{l = 0}^\infty \frac{t_2^l}{l!} (u_{b1}^*)^{k + l} \oint \dd w \, w^{k - l + r}
    = 2 \pi \ii (t_2 u_{b1}^*)^{1 + r} \sum_{k = 0}^{\infty} \frac{(t_1 t_2)^k}{k! \, (k + 1 + r)!},
\end{align*}
Note that we used $(u_{b1}^*)^2 = 1$.
The infinite series can be written as a Bessel-$I$ function, see \cite[\DLMFeqref{10.25}{2}]{NIST:DLMF}:
\begin{equation} \label{eq: w-sum Bessel I}
    \sum_{k = 0}^{\infty} \frac{x^k}{k! \, (k + 1 + r)!}
    = x^{-\frac{1+r}{2}} I_{1+r}\p*{2 \sqrt{x}}.
\end{equation}
In particular, we see that the sum over $u_{b1}^* = \pm 1$ reduces to the $u_{b1}^* = 1$ case times 2.
Using \cite[\DLMFeqref{10.29}{1}]{NIST:DLMF} we get $K_{-1}(z) = K_1(z)$, $I_{-1}(z) = I_1(z)$ and $K_{-2}(z) = \frac{2}{z} K_{-1}(z) + K_0(z)$.
The final result is
\begin{eqsplit}
    \frac{s_{12}}{s_{13} \sqrt{s_{22}}} \E\s*{\Tr\p*{\frac{1}{H - \eta \frac{s_{13} \sqrt{s_{22}}}{s_{12}} \zeta}}}
    \sim 2 \ii N^2 \Big[
        & K_0\p*{2 N \eta^{1/2} \sqrt{-\ii \zeta}} I_0\p*{2 N \eta^{1/2} \sqrt{-\ii \zeta}} \\
        & + K_1\p*{2 N \eta^{1/2} \sqrt{-\ii \zeta}} I_1\p*{2 N \eta^{1/2} \sqrt{-\ii \zeta}}
    \Big].
\end{eqsplit}

\section{Microscopic Scaling Limit of the One-Point Function}
\label{sec: Microscopic Scaling Limit}
In this section we prove Corollary~\ref{thm: microscopic one-point function}.
Integrating the asymptotic form of $\rho_{\infty}$ from \eqref{eq: density limit constant} yields
\[
    \int_{-\eta_N / 2}^{\eta_N / 2} \rho_{\infty}(E) \, \dd E
    \sim \int_{-\eta_N / 2}^{\eta_N / 2} \vartheta \abs{E}^{-\frac{\ell-1}{\ell+1}} \, \dd E
    = \vartheta (\ell + 1) \p*{\frac{\eta_N}{2}}^{\frac{2}{\ell+1}}.
\]
The right-hand side equals $\frac{1}{K N}$ if $\eta_N$ is given as in \eqref{eq: spacing scale}.
To determine $\vartheta$, we use \eqref{eq: Asymptotic density from saddle points} and Taylor-expand the solution for $z \to 0$.
The details for the last step depend on $K$ and are given below.

\subsection{Case \texorpdfstring{$K = 1$}{K = 1}}
Using $z = \ii \eta_N$ we see that the solution of $\nabla I(a) = 0$ satisfies
\[
    \frac{1}{a} = 1 \times a + \eta_N
    \qquad \iff \qquad
    0 = a^2 + \eta_N a - 1
\]
and the solution with $\Re(a(\eta_N)) > 0$ for $\eta_N > 0$ is $a(\eta_N) = \frac{1}{2} \p*{\sqrt{\eta_N^2 + 4} - \eta_N}$.
For $N \to \infty$ we have $\eta_N \to 0$ and $a(\eta_N) \to 1$, thus from \eqref{eq: Asymptotic density from saddle points} we get $\rho_{\infty}(0) = \frac{1}{\pi}$.
Since there is no divergence at 0, we have $\vartheta = \frac{1}{\pi}$ and $l = 1$, therefore $\eta_N = \frac{\pi}{N}$.

\subsection{Case \texorpdfstring{$K = 2$}{K = 2}}
From \eqref{eq: K=2 saddle-point equation} with $a_1 = \eta_N^{1/3} \wt{a}_1$, $a_2 = \eta_N^{-1/3} \wt{a}_2$ and $z = \ii \eta_N$ we get
\[
    0 = s_{11} \wt{a}_1^2 - \frac{1}{s_{12} \wt{a}_1 +\eta_N^{2/3}} + \eta_N^{2/3} \wt{a}_1, \qquad
    \wt{a}_2 = \frac{1}{s_{12} \wt{a}_1 +\eta_N^{2/3}}.
\]
We Taylor-expand the fractions in $\eta_N$ to see that $s_{11} \wt{a}_1^3 = \frac{1}{s_{12}} + O(\eta_N^{2/3})$.
Thus, $\wt{a}_1 = s_{11}^{-1/3} s_{12}^{-1/3} + O(\eta_N^{2/3})$ and $\wt{a}_2 = s_{11}^{1/3} s_{12}^{-2/3} + O(\eta_N^{2/3})$.
Extending analytically to the real line yields
\[
    a_2(E) = s_{11}^{1/3} s_{12}^{-2/3} (-\ii E)^{-1/3} + O(\abs{E}^{1/3}), \qquad
    \Re a_2(E) = \frac{\sqrt{3}}{2} s_{11}^{1/3} s_{12}^{-2/3} \abs{E}^{-1/3} + O(\abs{E}^{1/3}).
\]
The corresponding asymptotic density of states is then
\[
    \rho_{\infty}(E) = \frac{\sqrt{3}}{4 \pi} s_{11}^{1/3} s_{12}^{-2/3} \abs{E}^{-1/3} + O(\abs{E}^{1/3}),
\]
and therefore
\begin{equation} \label{eq: K=2 spacing scale}
    l = 2, \qquad
    \vartheta = \frac{\sqrt{3}}{4 \pi} s_{11}^{1/3} s_{12}^{-2/3}, \qquad
    \eta_N = 2 \frac{s_{12}}{\sqrt{s_{11}}} \bigg(\frac{2 \pi}{3\sqrt{3} } \bigg)^{3/2}N^{-3/2}.
\end{equation}

Note that the Meijer $G$-function has a branch cut along $\zeta \in \R$.
However, the real part of the $G$-function is continuous across this cut, so we can take the $\Im(\zeta) \downarrow 0$ limit of \eqref{eq: K=2 expectation limit} without trouble.

In Corollary~\ref{thm: microscopic one-point function} we state the asymptotic behavior for $\abs{\xi} \to 0$.
This follows from
\begin{align*}
    \MeijerG{3, 0}{0, 3}{-}{0, \frac{1}{2}, \frac{1}{2}}{x} & = \pi + 2 \sqrt{\pi x} \, \log(x) + O(\sqrt{x}), &
    \HypergeometricPFQ{0}{2}{-}{\frac{1}{2}, 1}{x} & = 1 + O(x), \\
    \MeijerG{3, 0}{0, 3}{-}{\frac{1}{2}, 1, \frac{3}{2}}{x} & = O(x), &
    \HypergeometricPFQ{0}{2}{-}{\frac{3}{2}, 2}{x} & = 1 + O(x), \\
    \MeijerG{3, 0}{0, 3}{-}{\frac{1}{2}, \frac{1}{2}, 1}{x} & = -\sqrt{\pi x} \, \log(x) + O(x), &
    \HypergeometricPFQ{0}{2}{-}{1, \frac{3}{2}}{x} & = 1 + O(x)
\end{align*}
as $x \to 0$.
The Meijer $G$ asymptotics can be derived from the contour integral representation \cite[\DLMFeqref{16.17}{1}]{NIST:DLMF} and the residue theorem.
The asymptotics of the generalized hypergeometric function follows directly from its series representation \cite[\DLMFeqref{16.2}{1}]{NIST:DLMF}.
The asymptotic behaviour for $\abs{\xi} \to \infty$ follows from
\begin{align*}
    \MeijerG{3, 0}{0, 3}{-}{0, \frac{1}{2}, \frac{1}{2}}{x} & = \frac{2 \pi}{\sqrt{3}} \ee^{-3 x^{1/3}} (1 + O(x^{-1/3})), &
    \HypergeometricPFQ{0}{2}{-}{\frac{1}{2}, 1}{x} & = \frac{1}{2 \sqrt{3 \pi}} x^{-1/6} \ee^{3 x^{1/3}} (1 + O(x^{-1/3})), \\
    \MeijerG{3, 0}{0, 3}{-}{\frac{1}{2}, 1, \frac{3}{2}}{x} & = \frac{2 \pi}{\sqrt{3}} x^{2/3} \ee^{-3 x^{1/3}} (1 + O(x^{-1/3})), &
    \HypergeometricPFQ{0}{2}{-}{\frac{3}{2}, 2}{x} & = \frac{1}{4 \sqrt{3 \pi}} x^{-5/6} \ee^{3 x^{1/3}} (1 + O(x^{-1/3})), \\
    \MeijerG{3, 0}{0, 3}{-}{\frac{1}{2}, \frac{1}{2}, 1}{x} & = \frac{2 \pi}{\sqrt{3}} x^{1/3} \ee^{-3 x^{1/3}} (1 + O(x^{-1/3})), &
    \HypergeometricPFQ{0}{2}{-}{1, \frac{3}{2}}{x} & = \frac{1}{4 \sqrt{3 \pi}} x^{-1/2} \ee^{3 x^{1/3}} (1 + O(x^{-1/3})).
\end{align*}
We refer to \cite[Chap.~V]{luke1969SpecialFunctionsApproximations} and \cite[\DLMFsecref{16.11}]{NIST:DLMF} for their source.

\subsection{Case \texorpdfstring{$K = 3$}{K = 3}}
Rescaling \eqref{eq: K=3 saddle-point equation} with $a_1 = \eta_N^{1/2} \wt{a}_1$, $a_2 = \wt{a}_2$, $a_3 = \eta_N^{-1/2} \wt{a}_3$ yields
\[
    \frac{1}{\wt{a}_1} = s_{11} \eta_N \wt{a}_1 + s_{12} \eta_N^{1/2} \wt{a}_2 + s_{13} \wt{a}_3 + \eta_N^{3/2}, \quad
    \frac{1}{\wt{a}_2} = s_{12} \eta_N^{1/2} \wt{a}_1 + s_{22} \wt{a}_2 + \eta_N, \quad
    \frac{1}{\wt{a}_3} = s_{13} \wt{a}_1 + \eta_N^{1/2}.
\]
We conclude $\wt{a}_2 = s_{22}^{-1/2} + O(\eta_N^{1/2})$ and from
\[
    0 = s_{11} \eta_N^{1/2} \wt{a}_1 + s_{12} \wt{a}_2 - \frac{1}{\wt{a}_1 (s_{13} \wt{a}_1 + \eta_N^{1/2})} + \eta_N
\]
that $\wt{a}_1 = (s_{12} s_{13} \wt{a}_2)^{-1/2} + O(\eta_N^{1/2}) = \frac{\sqrt[4]{s_{22}}}{\sqrt{s_{12} s_{13}}} + O(\eta_N^{1/2})$.
With $\wt{a}_3 = \frac{\sqrt{s_{12}}}{\sqrt[4]{s_{22}} \sqrt{s_{13}}} + O(\eta_N^{1/2})$ we find
\[
    \rho_{\infty}(E)
    = \frac{\sqrt{s_{12}}}{3 \pi \sqrt[4]{s_{22}} \sqrt{s_{13}}} \Re (-\ii E)^{-1/2} + O(1)
    = \frac{\sqrt{s_{12}}}{3 \pi \sqrt[4]{s_{22}} \sqrt{2 s_{13}}} \abs{E}^{-1/2} + O(1).
\]
The corresponding local spacing scale is
\begin{equation} \label{eq: K=3 spacing scale}
    l = 3, \qquad
    \vartheta = \frac{\sqrt{s_{12}}}{3 \pi \sqrt[4]{s_{22}} \sqrt{2 s_{13}}}, \qquad
    \eta_N = \frac{\pi^2}{4} \frac{s_{13} \sqrt{s_{22}}}{s_{12}} N^{-2}.
\end{equation}

In Corollary~\ref{thm: microscopic one-point function} we state the asymptotic behaviour for $\abs{\xi} \to 0$.
This follows from \cite[\DLMFeqref{10.31}{1} and \DLMFeqref{10.25}{1}]{NIST:DLMF}
\begin{align*}
    K_0(x) & = -\log(x) + \log(2) - \gamma + O(x^2 \log(x)), &
    I_0(x) & = 1 + O(x^2), \\
    K_1(x) & = \frac{1}{x} + \frac{x}{2} \log(x) + O(x), &
    I_1(x) & = \frac{x}{2} + O(x^2)
\end{align*}
as $x \to 0$, where $\gamma$ is the Euler–Mascheroni constant.
The asymptotics for $\abs{\xi} \to \infty$ follow from
\begin{align*}
    K_0(x) & = \sqrt{\frac{\pi}{2}} \, x^{-1/2} \ee^{-x} (1 + O(x^{-1})), &
    I_0(x) & = \frac{1}{\sqrt{2 \pi}} x^{-1/2} \ee^{x} (1 + O(x^{-1})), \\
    K_1(x) & = \sqrt{\frac{\pi}{2}} \, x^{-1/2} \ee^{-x} (1 + O(x^{-1})), &
    I_1(x) & = \frac{1}{\sqrt{2 \pi}} x^{-1/2} \ee^{x} (1 + O(x^{-1}))
\end{align*}
and we refer to \cite[\DLMFeqref{10.40}{1} and \DLMFeqref{10.40}{2}]{NIST:DLMF} for their source.

\section{Weak Non-Chirality Limit}
\label{sec: Weak Non-Chirality Limit}
In this section we prove Theorem~\ref{thm: microscopic resolvent expectation weak non-chirality}, i.e.\ we let $s_{11} \to 0$ (for $K = 2$) or $s_{12} \to 0$ (for $K = 3$) as $N \to \infty$ and derive the leading order of $\E\s*{\Tr((H - z)^{-1})}$ as $z \to 0$ with $N \to \infty$.
Note that the convergence rate of $z$ depends on the rate of $s_{11}$ and $s_{12}$, respectively.

\subsection{Case \texorpdfstring{$K = 2$}{K = 2}}
\label{sec: Weak Non-Chirality Limit for K = 2}
We start with the assumption that the convergence rates of $s_{11}$ and $z$ are $s_{11} \asymp N^{-\alpha}$ and $z \asymp N^{-\beta}$ for $\alpha, \beta > 0$.
Later we will choose the rate $\beta$ depending on $\alpha$ in such a way that we get a non-trivial scaling limit.
From Section~\ref{sec: Large N Limit at the Origin} we already know that $\beta = \frac{3}{2}$ is the appropriate rate when $\alpha = 0$.
Furthermore, by comparing with the chiral GUE we see that $\beta = 1$ for $\alpha \to \infty$.
Let the block variance profile be
\[
    S = \begin{pmatrix}
        N^{-\alpha} \sigma_{11} & 1 \\
        1 & 0
    \end{pmatrix}, \qquad
    \sigma_{11} \in [0, \infty),\ \alpha \geq 0
\]
and $z = N^{-\beta} \zeta$ with $\Im(\zeta) > 0$ and $\beta \geq 1$.

Analogous to Section~\ref{sec: Large N Limit at the Origin}, we get the saddle-point equations
\[
    \frac{1}{a_1} = -\ii N^{-\beta} \zeta + N^{-\alpha} \sigma_{11} + a_2, \qquad
    \frac{1}{a_2} = -\ii N^{-\beta} \zeta + a_1.
\]
Further assuming that $a_1 \asymp N^{-c}$ for $c < \beta$ leads to $a_2 \asymp N^c$ via the second equation.
By comparing the first equation times $a_1$ and the second equation times $a_2$ we obtain
\[
    a_1
    = \frac{-\ii N^{-\beta} \zeta a_2}{-\ii N^{-\beta} \zeta + N^{-\alpha} \sigma_{11} a_1}
    \asymp \frac{N^{-\beta + c}}{N^{-\beta} + N^{-\alpha - c}}
    \sim N^{-\beta + c - \max(-\beta, -\alpha - c)}.
\]
Since we started with $a_1 \asymp N^{-c}$, we choose $c = \max(-c, \beta - 2 c - \alpha)$ to get a non-trivial solution for the saddle-point equations.
Further simplification leads to $c = \frac{1}{3} \max(\beta - \alpha, 0)$.

As in section~\ref{sec: Large N Limit at the Origin} we transform the integration variables $(a_1, a_2) = (N^{-c} u_a v, N^c \frac{u_a}{v})$ for $u_a, v > 0$ and $(b_1, b_2) = (N^{-c} u_b w, N^c \frac{u_b}{w})$ for $u_b, w \in U(1)$ and arrive at
\begin{eqsplit*}
    \E\s*{\Tr\p*{\frac{1}{H - N^{-\beta} \zeta}}}
    & = - \frac{2 \ii N^3}{(2 \pi)^2} \int_0^\infty \dd v \oint \dd w
    \int_0^\infty \dd u_a \oint \dd u_b \, \ee^{-N J(u)} F(u; v, w) D(u; v, w),
\end{eqsplit*}
where
\begin{align*}
    J(u) & := u_a^2 - 2 \log(u_a) - u_b^2 + 2 \log(u_b),
    \\
    L(u; v, w) & := N^{1 - \alpha - 2 c} \frac{\sigma_{11}}{2} (u_a^2 v^2 - u_b^2 w^2) + N^{1 - \beta + c} (-\ii \zeta) \p*{\frac{u_a}{v} - \frac{u_b}{w}} + N^{1 - \beta - c} (-\ii \zeta) (u_a v - u_b w),
    \\
    F(u; v, w) & := N^c \ee^{-L(u; v, w)} \frac{1 + N^{-2 c} v^2}{v^2 w} u_a^2 u_b,
    \\
    D(u_a, u_b; v, w) & := \det\s*{ S + \diag\p*{\frac{N^{2 c}}{u_a v u_b w}, \frac{N^{-2 c} v w}{u_a u_b}} }
    = D^{(0)}(u) + N^{-\alpha - 2 c} D^{(1)}(u; v, w)
\end{align*}
and
\[
    D^{(0)}(u) := \frac{1}{u_a^2 u_b^2} - 1, \qquad
    D^{(1)}(u; v, w) := \sigma_{11} \frac{v w}{u_a u_b}.
\]
To get a non-trivial scaling limit, at least one of the powers of $N$ in the $L$-function must be zero.
Setting $1 - \alpha - 2 c = 0$ and $1 - \beta + c = 0$ with $c = \frac{1}{3} \max(\beta - \alpha, 0)$ gives the solution
\[
    \begin{cases}
        \beta = \frac{1}{2} (3 - \alpha),\; c = \frac{1}{3} (\beta - \alpha) & \text{for } \alpha < 1, \\
        \beta = 1, \; c = 0 & \text{for } \alpha \geq 1.
    \end{cases}
\]
In particular, when $\alpha > 1$ the term in $L(u; v, w)$ which is proportional to $N^{1 - \alpha - 2 c}$ will be subleading and the integral over $N^{-\alpha - 2 c} D^{(1)}(u; v, w)$ will not contribute to the leading order.

We perform a saddle-point approximation to obtain
\begin{eqsplit*}
    \int_0^\infty \dd u_a \oint \dd u_b \, \ee^{-N J} F D^{(0)}
    & \sim \frac{\ii \pi}{4 N^2} N^c \ee^{
        -\frac{1}{2} \sigma_{11} v^2 + \ii \zeta \p*{\frac{1}{v} + N^{-2 c} v}
    } \sum_{u_b^* = \pm 1} \ee^{
        \frac{1}{2} \sigma_{11} w^2 - \ii \zeta \p*{\frac{1}{w} + N^{-2 c} w} u_b^*
    }
    \\
    & \qquad \times \frac{1 + N^{-2 c} v^2}{v^2 w} \bigg[
        \sigma_{11} v^2 + (-\ii \zeta) \frac{1 + N^{-2 c} v^2}{v} - 1
    \\
    & \hphantom{\qquad \times \frac{1 + N^{-2 c} v^2}{v^2 w} \bigg[}
        + \sigma_{11} w^2 + (-\ii \zeta) \frac{1 + N^{-2 c} w^2}{w} u_b^*
    \bigg]
\end{eqsplit*}
and
\begin{eqsplit*}
    \int_0^\infty \dd u_a \oint \dd u_b \, \ee^{-N J} F D^{(1)}
    & \sim \frac{\ii \pi}{2 N} N^c \ee^{
        -\frac{1}{2} \sigma_{11} v^2 + \ii \zeta \p*{\frac{1}{v} + N^{-2 c} v}
    } \sum_{u_b^* = \pm 1} \ee^{
        \frac{1}{2} \sigma_{11} w^2 - \ii \zeta \p*{\frac{1}{w} + N^{-2 c} w} u_b^*
    }
    \\
    & \qquad \times \frac{1 + N^{-2 c} v^2}{v} \sigma_{11} u_b^*.
\end{eqsplit*}
Next, we want to evaluate the $\int_0^{\infty} \dd v$ and $\oint \dd w$ integrals.
After we expand the integrand, we see that the $v$-integrals are of the form
\[
    \int_0^{\infty} \dd v \, \ee^{-t_1 v^2 + t_2 v + t_3 \frac{1}{v}} v^n \qquad \text{with } t_1 = \frac{1}{2} \sigma_{11},\ t_2 = \ii \zeta N^{-2 c},\ t_3 = \ii \zeta \text{ and } n \in \Z.
\]
If $t_1 = 0$ (i.e.\ when $\sigma_{11} = 0$), we get a result in terms of Bessel-$K$ functions, see \eqref{eq: v-integral Bessel K}.
If $t_2 = 0$ (i.e.\ when $c > 0$ and $N \to \infty$), we get a result in terms of Meijer $G$-functions, see \eqref{eq: v-integral Meijer G}.
But for the general case, we did not find a formula in terms of special functions.
Thus, we leave this integral unevaluated in Theorem~\ref{thm: microscopic resolvent expectation weak non-chirality}.
In contrast, the $w$-integral can be performed as in Section~\ref{sec: Large N Limit at the Origin}
\begin{align*}
    \int_{U(1)} \dd w \, \ee^{t_1 w^2 + t_2 w + t_3 \frac{1}{w}} w^r
    = 2 \pi \ii \sum_{k, l = 0}^{\infty} \frac{t_1^k t_2^l t_3^{2 k + l + 1 + r}}{k! \, l! \, (2 k + l + 1 + r)!}, \quad r \in \Z.
\end{align*}
Using this formula, it is easy to see that the sum over $u_b^* = \pm 1$ simplifies to $2$ times the $u_b^* = 1$ case.
If $t_1 = 0$, only the $k = 0$ terms are non-zero and the sum over $l$ gives a Bessel-$I$ function, see \eqref{eq: w-sum Bessel I}.
If $t_2 = 0$, only the $l = 0$ terms are non-zero and the sum over $k$ can be written in terms of a generalized hypergeometric function, see \eqref{eq: w-sum Hypergeometric F}.
However, we did not find a formula to express the general double-sum in terms of special functions.

For the possible cases of $\alpha$ we arrive at the following results:
\begin{description}
    \item[$\alpha < 1$]: Then $\beta > \alpha$ and $c > 0$ and the large-$N$ limit of the expectation of the resolvent is the same as in Theorem~\ref{thm: microscopic resolvent expectation} with an additional $N^{-\alpha/2}$ factor.

    \item[$\alpha = 1$]: Then $\beta = 1$ and $c = 0$ and we get the formula in Theorem~\ref{thm: microscopic resolvent expectation weak non-chirality}.
    The coefficients $c_{n, r}$ follow from the expansion
    \[
        \frac{1 + v^2}{v^2 w} \s*{
            \sigma_{11} v^2 + (-\ii \zeta) \frac{1 + v^2}{v} - 1
            + \sigma_{11} w^2 + (-\ii \zeta) \frac{1 + w^2}{w}
            + 2 \sigma_{11} v w
        }
        = \sum_{n, r \in \Z} c_{n, r} v^n w^r.
    \]
    In the special case $\sigma_{11} = 0$ we have
    \[
        \E\s*{\Tr\p*{\frac{1}{H - N^{-1} \zeta}}}
        \sim \frac{1}{4 \pi} N \int_0^{\infty} \dd v \, \ee^{\ii \zeta (\frac{1}{v} + v)}
        \oint \dd w \, \ee^{\ii \zeta (\frac{1}{w} + w)}
        \frac{1 + v^2}{v^2 w} \s*{(-\ii \zeta) \p2{\frac{1}{v} + v + \frac{1}{w} + w} - 1}.
    \]
    The $v$- and $w$-integral can be evaluated via \eqref{eq: v-integral Bessel K} and \eqref{eq: w-sum Bessel I}, respectively.
    Finally, we use \cite[\DLMFeqref{10.29}{1}]{NIST:DLMF} to simplify the resulting combinations of Bessel-$K$ and Bessel-$I$ functions and arrive at \eqref{eq: chiralGUE expectation of resolvent}.

    \item[$\alpha > 1$]: Then $\alpha > \beta$ and $c = 0$ and thus for $N \to \infty$ we get the same large-$N$ limit as for $\alpha = 1$ and $\sigma_{11} = 0$.
\end{description}

\subsection{Case \texorpdfstring{$K = 3$}{K = 3}}
\label{sec: Weak Non-Chirality Limit for K = 3}
Let
\[
    S = \begin{pmatrix}
        0 & N^{-\alpha} \sigma_{12} & 1 \\
        N^{-\alpha} \sigma_{12} & 1 & 0 \\
        1 & 0 & 0
    \end{pmatrix}, \qquad
    \sigma_{12} \in [0, \infty),\ \alpha \geq 0
\]
and $z = N^{-\beta} \zeta$ with $\Im(\zeta) > 0$ and $\beta \geq 1$.
Then the saddle-point equations are
\[
    \frac{1}{a_1} = -\ii N^{-\beta} \zeta + N^{-\alpha} \sigma_{12} a_2 + a_3, \qquad
    \frac{1}{a_2} = -\ii N^{-\beta} \zeta + N^{-\alpha} \sigma_{12} a_1 + a_2, \qquad
    \frac{1}{a_3} = -\ii N^{-\beta} \zeta + a_1.
\]
Assuming $a_1 \asymp N^{-c}$ with $c < \beta$ we get $a_3 \asymp N^c$ from the third equation.
The second equation can only be satisfied if $a_2 \asymp 1$ (assuming $c \geq -\alpha$).
By comparing the first equation times $a_1$ and the third equation times $a_3$ we get $c = \frac{1}{2} \max(\beta - \alpha, 0)$.
This is also consistent with $\beta = 2$ for $\alpha = 0$ (from section~\ref{sec: Large N Limit at the Origin}) and $\beta = 1$ for $\alpha \to \infty$ (the limiting density is the semi-circle law).

We transform the integration variables
\begin{eqsplit*}
    (a_1, a_2, a_3) & = \p*{N^{-c} u_{a1} v, u_{a2}, N^c \frac{u_{a1}}{v}} \qquad \text{for } u_{a1}, u_{a2}, v > 0,
    \\
    (b_1, b_2, b_3) & = \p*{N^{-c} u_{b1} w, u_{b2}, N^c \frac{u_{b1}}{w}} \qquad \text{for } u_{b1}, u_{b2}, w \in U(1)
\end{eqsplit*}
and arrive at
\begin{eqsplit*}
    \E\s*{\Tr\p*{\frac{1}{H - N^{-\beta} \zeta}}}
    & = -\frac{N^4}{4 \pi^3}
        \int_0^\infty \dd v \int_{U(1)} \dd w
        \int_{(0, \infty)^2} \dd u_a \oint^2 \dd u_b \,
        \ee^{-N J(u)} F(u; v, w) D(u; v, w),
\end{eqsplit*}
where
\begin{align*}
    J(u) & := u_{a1}^2 - 2 \log(u_{a1}) + \frac{1}{2} u_{a2}^2 - \log(u_{a2}) - u_{b1}^2 + 2 \log(u_{b1}) - \frac{1}{2} u_{b2}^2 + \log(u_{b2}),
    \\
    L(u; v, w) & := N^{1 - \alpha - c} \sigma_{12} (u_{a1} u_{a2} v - u_{b1} u_{b2} w)
    + N ^{1 - \beta + c} (-\ii \zeta) \p2{\frac{u_{a1}}{v} - \frac{u_{b1}}{w}}
    \\
    & \quad + N ^{1 - \beta} (-\ii \zeta) (u_{a2} - u_{b2})
    + N^{1 - \beta - c} (-\ii \zeta) (u_{a1} v - u_{b1} w),
    \\
    F(u; v, w) & := \ee^{-L(u; v, w)} \frac{N^{-c} u_{a1} v^2 + u_{a2} v + N^c u_{a1}}{v^2 w} u_{a1} u_{b1},
    \\
    D(u; v, w) & := \det\s*{ S + \diag\p*{\frac{N^{2 c}}{u_{a1} v u_{b1} w}, \frac{1}{u_{a2} u_{b2}}, \frac{N^{-2 c} v w}{u_{a1} u_{b1}}} }
    = D^{(0)}(u) + N^{-2 \alpha - 2 c} D^{(1)}(u; v, w).
\end{align*}
and
\[
    D^{(0)}(u) := \p*{\frac{1}{u_{a1}^2 u_{b1}^2} - 1} \p*{\frac{1}{u_{a2} u_{b2}} + 1}, \qquad
    D^{(1)}(u; v, w) := -\sigma_{12}^2 \frac{v w}{u_{a1} u_{b1}}.
\]
Examining the powers of $N$ in the $L$-function, we see that $\beta = 2 - \alpha$, $c = 1 - \alpha$ are the appropriate rates when $\alpha < 1$ and $\beta = 1$, $c = 0$ are the rates for $\alpha \geq 1$.

A saddle-point approximation yields
\begin{eqsplit*}
    \int \dd u \, \ee^{-N J} F D^{(0)}
    & \sim - \frac{\pi^2}{2 N^3} N^c \sum_{u_{b1}^* = \pm 1} \ee^{-L(1, 1, u_{b1}^*, 1; v, w)} \\
    & \quad \times \Bigg[ \p*{
        N^{1 - \alpha - c} \sigma_{12} v
        + N^{1 - \beta + c} (-\ii \zeta) \frac{1}{v}
        + N^{1 - \beta - c} (-\ii \zeta) v
    } \frac{1 + N^{-c} v + N^{-2 c} v^2}{v^2 w} \\
    & \qquad + \p*{
        N^{1 - \alpha - c} \sigma_{12} w
        + N^{1 - \beta + c} (-\ii \zeta) \frac{1}{w}
        + N^{1 - \beta - c} (-\ii \zeta) w
    } \frac{1 + N^{-c} v + N^{-2 c} v^2}{v^2 w} u_{b1}^* \\
    & \qquad - \frac{1 + N^{-2 c} v^2}{w}
    \Bigg].
\end{eqsplit*}
Observe that $\int \dd u \, \ee^{-N J} F D^{(1)} = O(N^{-2 + c})$, therefore the $N^{-2 \alpha - 2 c} D^{(1)}$ term cannot contribute to the leading order of the final result.

Next, we evaluate the $\int_0^{\infty} \dd v$ and $\oint \dd w$ integrals with \eqref{eq: v-integral Bessel K} and \eqref{eq: w-sum Bessel I}.
Using \cite[\DLMFeqref{10.29}{1}]{NIST:DLMF}, we simplify the resulting Bessel-$K$ and Bessel-$I$ functions and arrive at
\begin{eqsplit*}
    \frac{1}{N} \E\s*{\Tr\p*{\frac{1}{H - N^{-\beta} \zeta}}}
    & \sim 2 \ii N^2 \Big[
        (N^{-\alpha} \sigma_{12} - 2 N^{-\beta} \ii \zeta) \p1{K_0(x) I_0(x) + K_1(x) I_1(x)}
    \\
    & \hphantom{= 2 \ii N^2 \Big[}
        + \sqrt{N^{-\beta} (-\ii \zeta) (N^{-\alpha} \sigma_{12} - N^{-\beta} \ii \zeta)} \p1{K_1(x) I_0(x) + K_0(x) I_1(x)}
    \Big]
\end{eqsplit*}
where $x = 2 N \sqrt{N^{-\beta} (-\ii \zeta) (N^{-\alpha} \sigma_{12} - N^{-\beta} \ii \zeta)}$.
Using the formula \cite[\DLMFeqref{10.28}{2}]{NIST:DLMF} for the Wronskian of $I_{\nu}$ and $K_{\nu}$ we see that the term on the second line is equal to $\frac{1}{2 N}$.

If $\alpha < 1$, then $\beta = 2 - \alpha$ and $c = 1 - \alpha$ and for $N \to \infty$ we arrive again at the result for the strong non-chirality limit from Theorem~\ref{thm: microscopic resolvent expectation}.
If $\alpha = 1$, then $\beta = 1$ and $c = 0$ and we obtain the formula for Theorem~\ref{thm: microscopic resolvent expectation weak non-chirality}.
The case $\alpha > 1$ yields the same large-$N$ limit as $\alpha = 1$ with $\sigma_{12} = 0$.

\begin{appendices}

\section{Multi-Dimensional Saddle-Point Method}
\label{app: Saddle-Point Method}
Let us explain the one-dimensional saddle-point method first.
Consider the integral $I_N := \int_{\gamma} \dd u \, \ee^{N S(u)} f(u)$,
where $\gamma$ is a contour in $\C$ and $S, f: \C \to \C$ are analytic.
We want to find an asymptotic expansion of $I_N$ as $N \to \infty$.
The idea of the \emph{method of steepest descent} (also known as the \emph{saddle-point method}) is to deform the contour $\gamma$ into a contour $\tilde{\gamma}$ without changing the endpoints and such that
\begin{enumerate}[label=(\roman*)]
    \item $\tilde{\gamma}$ passes through a critical point $u^*$ of $S$ (i.e.\ $S'(u^*) = 0$),

    \item $\Im(S(u)) = \Im(S(u^*))$ holds on some neighbourhood of $u^*$ on $\tilde{\gamma}$,

    \item $\Re(S(u)) < \Re(S(u^*))$ for $u \neq u^*$ on $\tilde{\gamma}$ (or equivalently, $\max_{u \in \tilde{\gamma}} \abs{\ee^{N S(u)}}$ is attained only at $u^*$).
\end{enumerate}
Let $\tilde{\gamma}_{\eps}$ be a portion of $\tilde{\gamma}$ in a neighbourhood of $u^*$ where (ii) holds.
Then $I_N = I_{N, \eps} + I_{N, \eps}'$, where $I_{N, \eps} := \int_{\tilde{\gamma}_{\eps}} \dd u \, \ee^{N S(u)} f(u)$ and $I_{N, \eps}' := \int_{\tilde{\gamma} \setminus \tilde{\gamma}_{\eps}} \dd u \, \ee^{N S(u)} f(u)$.
We can estimate $\abs{I_{N, \eps}'} \leq C \ee^{-N \delta} \abs{\ee^{N S(u^*)}}$ for some constants $C > 0$ and $\delta > 0$ (for example, $\delta = - \sup_{u \in \tilde{\gamma} \setminus \tilde{\gamma}_{\eps}} \Re(S(u) - S(u^*))$).
For $I_{N, \eps}$ we first parametrize $\tilde{\gamma}_{\eps}$ by $u(x) = u^* + \mu x$, where $x \in (-\eps, \eps)$ and $\mu \in U(1)$ is the contour direction.
We extend the integration range from $x \in (-\eps, \eps)$ to $x \in (-\infty, \infty)$, this introduces only an error of size $\abs{\ee^{N S(u^*)}} O(\ee^{-N \eps})$.
Finally, we use Laplace's method the derive the asymptotic expansion
\[
    I_N = \sqrt{\frac{2 \pi}{N}} \frac{1}{\sqrt{-S''(u^*)}} \ee^{N S(u^*)} \p*{f(u^*) + O\p*{N^{-1}}},
    \qquad \text{where } \arg \sqrt{-S''(u^*)} \in \p*{-\frac{\pi}{4}, \frac{\pi}{4}}.
\]
More generally, if there are several discrete critical points $u^*$ where $\Re S$ attains its maximum and that can be reached by deforming the integration contour, then each of these points contributes to the final result.
A similar asymptotic formula exists when the maximum of $\Re S$ is attained at the one of the endpoints of $\gamma$.

For the integrals in Section~\ref{sec: Large N Limit at the Origin} we need the $d$-dimensional version of the saddle-point method including the lower order terms.
The asymptotic expansion in this case has been derived before by Fedoryuk \cite{fedoryuk1977saddle} (see also \cite{fedoryuk1989asymptotic}) and reads
\[
    \int_{\gamma} \dd u \, \ee^{N S(u)} f(u)
    = \p*{\frac{2 \pi}{N}}^{d/2} \frac{1}{\sqrt{\det(-S''(u^*))}} \ee^{N S(u^*)} \p*{f(u^*) + \sum_{k = 1}^\infty c_k N^{-k}},
\]
where $c_k$ depends on $S, f$ and their derivatives at $u^*$.
The choice of branch for $\sqrt{\det(-S''(u^*))}$ depends on the orientation of the contour.

In the following we present a derivation of the asymptotic expansion for the integral $\int_{\gamma} \dd u \, \ee^{-N J(u)} F(u) D(u)$, where $\gamma$ is a $d$-dimensional contour and $J, F, D: \C^d \to \C$ are analytic in the region where we will deform the contour.
Let $J$ have (at least) one non-degenerate critical point $u^*$ that is not at the boundary of $\gamma$ and let $F(u^*) \neq 0$ while $D$ and derivatives of $D$ can be zero at $u^*$.
Furthermore, we assume that $\partial_i \partial_j J(u^*) = 0$ for $i \neq j$ (for example, this is the case when $J(u) = \wt{J}_1(u_1) + \dots + \wt{J}_d(u_d)$).
This last assumption simplifies the parametrization of $\tilde{\gamma}_{\eps}$ and helps selecting the correct branch of $\sqrt{\det J''(u^*)}$ in the final expression.

\paragraph{Step 1: Critical points of \texorpdfstring{$J$}{J}}
Let $u^* \in \C^d$ be a non-degenerate critical point of $J$, i.e.\ $\partial_i J(u^*) = 0$ for all $i = 1, \dots, d$ and $\det J''(u^*) \neq 0$, where $\partial_i$ is the derivative with respect to the variable $u_i$ and $J'' = (\partial_i \partial_j J)_{i,j}$ is the Hessian.
We deform the contour $\gamma$ into $\tilde{\gamma}$ so that (i), (ii) and (iii) from above hold and denote by $\tilde{\gamma}_{\eps}$ a neighbourhood of $u^*$ in $\tilde{\gamma}$ where (ii) holds.
Let $\mu \in \C^d$ with $\abs{\mu_i} = 1$ and such that $\mu_i^2 \, \partial_i^2 J(u^*) > 0$.
In practice we use $\mu_i = \pm \ee^{-\frac{\ii}{2} \arg(\partial_i^2 J(u^*))}$ and choose the sign such that $\mu_i$ has the same orientation as the contour $\tilde{\gamma}_{\eps}$.
Now we parametrize $\tilde{\gamma}_{\eps}$ as $u = u^* + \frac{\mu y}{\sqrt{N}}$, where $y \in (-\sqrt{N} \, \eps, \sqrt{N} \, \eps)^d$ and $\mu y := (\mu_i y_i)_{i = 1, \dots, d}$ is a shorthand for entrywise multiplication of two vectors.
Then
\begin{eqsplit*}
    & \int_{\gamma} \dd u \, \ee^{-N J(u)} F(u) D(u)
    = \int_{\tilde{\gamma}_{\eps}} \dd u \, \ee^{-N J(u)} F(u) D(u) + \int_{\tilde{\gamma} \setminus \tilde{\gamma}_{\eps}} \dd u \, \ee^{-N J(u)} F(u) D(u)
    \\
    & \quad = \p*{\prod_{i = 1}^{d} \frac{\mu_i}{\sqrt{N}}} \p*{
        \int_{\R^d} \dd y \, \ee^{-N J\p*{u^* + \frac{\mu y}{\sqrt{N}}}} F\p*{u^* + \frac{\mu y}{\sqrt{N}}} D\p*{u^* + \frac{\mu y}{\sqrt{N}}}
    }
    \\
    & \qquad + O(N^{-d/2} \ee^{-N \eps - \Re(J(u^*))})
    + O(\ee^{-N \delta - \Re(J(u^*))})
\end{eqsplit*}
for some $\delta > 0$.

\paragraph{Step 2: Taylor expansion}
For a multi-index $\alpha \in \N_0^d$ define $\abs{\alpha} := \sum_{i = 1}^d \alpha_i$, $\alpha ! := \prod_{i = 1}^d \alpha_i !$, $(\mu y)^\alpha := \prod_{i = 1}^d (\mu_i y_i)^{\alpha_i}$ and $\partial^\alpha := \prod_{i = 1}^d \partial_i^{\alpha_i}$.
Then introduce for $k \in \N_0$ the polynomials
\[
    J_k(y) := \sum_{\abs{\alpha} = k} \frac{(\mu y)^\alpha}{\alpha !} \partial^\alpha J(u^*), \qquad
    F_k(y) := \sum_{\abs{\alpha} = k} \frac{(\mu y)^\alpha}{\alpha !} \partial^\alpha F(u^*), \qquad
    D_k(y) := \sum_{\abs{\alpha} = k} \frac{(\mu y)^\alpha}{\alpha !} \partial^\alpha D(u^*).
\]
Then $J\p*{u^* + \frac{\mu y}{\sqrt{N}}} = \sum_{k = 0}^{\infty} N^{-k/2} J_k(y)$ and similar for $F$ and $D$.
Via Taylor-expansion we obtain
\begin{eqsplit*}
    & \int_{\R^d} \dd y \, \ee^{-N J\p*{u^* + \frac{\mu y}{\sqrt{N}}}} F\p*{u^* + \frac{\mu y}{\sqrt{N}}} D\p*{u^* + \frac{\mu y}{\sqrt{N}}}
    \\
    & \quad = \ee^{-N J_0} \int_{\R^d} \dd y \, \ee^{-J_2} \bigg(
        F_0 D_0
        + \frac{1}{\sqrt{N}} \s1{
            F_1 D_0 + F_0 D_1 - F_0 D_0 J_3
        }
    \\
    & \hphantom{\quad = \ee^{-N J_0} \int_{\R^d} \dd y \, \ee^{-J_2} \bigg(}
        + \frac{1}{N} \s2{
                F_2 D_0 + F_1 D_1 + F_0 D_2 - F_1 D_0 J_3 - F_0 D_1 J_3 + \frac{1}{2} F_0 D_0 J_3^2 - F_0 D_0 J_4
            }
    \\
    & \hphantom{\quad = \ee^{-N J_0} \int_{\R^d} \dd y \, \ee^{-J_2} \bigg(}
        + O\p*{\frac{1}{N^{3/2}}}
    \bigg),
\end{eqsplit*}
where $J_k \equiv J_k(y)$, $F_k \equiv F_k(y)$ and $D_k \equiv D_k(y)$ for brevity.

\paragraph{Step 3: Integral over polynomials}
For $\Re(a) > 0$ the Gaussian integral over a monomial is
\[
    \int_{\R} \dd y \, \ee^{-\frac{1}{2} a y^2} y^k = \begin{cases}
        \sqrt{2 \pi} \, (k - 1)!! \, a^{-\frac{k + 1}{2}}, & k \text{ even}, \\
        0, & k \text{ odd}.
    \end{cases}
\]
If $J_2$ contains no mixed derivatives, i.e.\ $\partial_i \partial_j J(u^*) = 0$ for $i \neq j$, then
\[
    \int_{\R^d} \dd y \, \ee^{-J_2} y^\alpha = 0
    \qquad \text{if and only if some entry $\alpha_i$ is odd}.
\]
In particular, in the Taylor expansion any term of order $N^{-k/2}$ with $k$ odd will integrate to 0, e.g.
\[
    \int_{\R^d} \dd y \, \ee^{-J_2} \frac{1}{\sqrt{N}} \s1{F_1 D_0 + F_0 D_1 - F_0 D_0 J_3} = 0.
\]

\paragraph{Step 4: Compute the leading order}
Assume that $F_0 \neq 0$ and define $h_i := \partial_i^2 J(u^*)$ for brevity.

If $D_0 \neq 0$, then for $N \to \infty$ the leading order is
\begin{eqsplit*}
    & \int_{\R^d} \dd y \, \ee^{-N J\p*{u^* + \frac{\mu y}{\sqrt{N}}}} F\p*{u^* + \frac{\mu y}{\sqrt{N}}} D\p*{u^* + \frac{\mu y}{\sqrt{N}}}
    = \ee^{-N J_0} \p*{\int_{\R^d} \dd y \, \ee^{-J_2} F_0 D_0 + O\p*{\frac{1}{N}}}
    \\
    & \quad = (2 \pi)^{d/2} \ee^{-N J(u^*)} \p*{\prod_{i = 1}^{d} \abs{h_i}^{-1/2}} \p*{F(u^*) D(u^*) + O\p*{\frac{1}{N}}}.
\end{eqsplit*}

If $D_0 = 0$, then $F_2 D_0 - F_1 D_0 J_3 + \frac{1}{2} F_0 D_0 J_3^2 - F_0 D_0 J_4 = 0$ and thus
\begin{eqsplit*}
    & \int_{\R^d} \dd y \, \ee^{-N J\p*{u^* + \frac{\mu y}{\sqrt{N}}}} F\p*{u^* + \frac{\mu y}{\sqrt{N}}} D\p*{u^* + \frac{\mu y}{\sqrt{N}}}
    \\
    & \quad = \ee^{-N J_0} \p*{\int_{\R^d} \dd y \, \ee^{-J_2} \frac{1}{N} \s*{F_1 D_1 + F_0 D_2 - F_0 D_1 J_3} + O\p*{\frac{1}{N^2}}}
    \\
    & \quad = (2 \pi)^{d/2} \frac{1}{N} \ee^{-N J(u^*)} \p*{\prod_{i = 1}^{d} \abs{h_i}^{-1/2}} \\& \qquad \times \p*{\sum_{i = 1}^d \s2{
        \frac{1}{h_i} \partial_i F(u^*) \, \partial_i D(u^*)
        + \frac{1}{2 h_i} F(u^*) \, \partial_i^2 D(u^*)
        - \frac{1}{2 h_i^2} F(u^*) \, \partial_i D(u^*) \, \partial_i^3 J(u^*)
    } + O\p*{\frac{1}{N^2}}}.
\end{eqsplit*}

\paragraph{Summary}
As $N \to \infty$ the integral has an asymptotic expansion
\[
    \int_{\gamma} \dd u \, \ee^{-N J(u)} F(u) D(u)
    = \p*{\frac{2 \pi}{N}}^{d/2} \ee^{-N J(u^*)} \p*{\prod_{i = 1}^{d} \frac{\mu_i}{\sqrt{\abs{\partial_i^2 J(u^*)}}}} \times \p*{F(u^*) D(u^*) + \sum_{k = 1}^\infty c_k N^{-k}}.
\]
If $D(u^*) = 0$, then
\[
    c_1 = \sum_{i = 1}^d \s*{
        \frac{\partial_i F(u^*) \, \partial_i D(u^*)}{\partial_i^2 J(u^*)}
        + \frac{F(u^*) \, \partial_i^2 D(u^*)}{2 \, \partial_i^2 J(u^*)}
        - \frac{F(u^*) \, \partial_i D(u^*) \, \partial_i^3 J(u^*)}{2 (\partial_i^2 J(u^*))^2}
    }.
\]
The leading term is bounded by $O(\max_{u^*} \abs{\ee^{-N J(u^*)}} N^{-d/2 - \ceil{L/2}})$, where $L \equiv L(u^*) \in \N_0$ is such that $D_l = 0$ for $l < L$ and $D_L \neq 0$ at $u^*$.
The sub-leading terms are bounded by $O(\max_{u^*} \abs{\ee^{-N J(u^*)}} N^{-d/2 - (\ceil{L/2} + 1)})$.

\end{appendices}

\printbibliography

\end{document}